\documentclass[reqno,11pt]{amsart}
\usepackage{mathrsfs}
\usepackage{verbatim}
\usepackage{upgreek}
\usepackage{amsmath}
\usepackage{amsxtra}
\usepackage{amscd}
\usepackage{amsthm}
\usepackage{amsfonts}
\usepackage{amssymb}
\usepackage{eucal}
\usepackage{tikz-cd}
\usepackage[matrix,arrow,curve]{xy}
\usepackage{young}
\usepackage[vcentermath]{youngtab}
\usepackage{mathtools}

\theoremstyle{plain}
\newtheorem{Thm}{Theorem}[section]
\newtheorem{Cor}[Thm]{Corollary}
\newtheorem{Lem}[Thm]{Lemma}
\newtheorem{Prop}[Thm]{Proposition}
\newtheorem{Conj}[Thm]{Conjecture}
\newtheorem{Warning}[Thm]{Warning}
\newtheorem{Assumption}[Thm]{Assumption}

\newtheorem{thmx}{Theorem}

\theoremstyle{definition}
\newtheorem{Def}[Thm]{Definition}

\theoremstyle{definition}
\newtheorem{Ex}[Thm]{Example}

\theoremstyle{remark}
\newtheorem{Rem}[Thm]{Remark}

\newtheorem{conj}{Conjecture}

\numberwithin{equation}{section}
\renewcommand{\rm}{\normalshape}

\newif\ifShowLabels
\ShowLabelstrue
\newdimen\theight
\def\TeXref#1{%
	\leavevmode\vadjust{\setbox0=\hbox{{\tt
				\quad\quad  {\small \rm #1}}}%
		\theight=\ht0
		\advance\theight by \lineskip
		\kern -\theight \vbox to
		\theight{\rightline{\rlap{\box0}}%
			\vss}%
}}%

\ShowLabelsfalse

\newenvironment{thm}[1]%
{ \begin{Thm} \label{T:#1}  \ifShowLabels \TeXref{T:#1} \fi }%
	{ \end{Thm} }

\renewcommand{\th}[1]{\begin{thm}{#1} \sl }
	\renewcommand{\eth}{\end{thm} }

\newenvironment{lemma}[1]%
{ \begin{Lem} \label{L:#1}  \ifShowLabels \TeXref{L:#1} \fi }%
	{ \end{Lem} }
\newcommand{\lem}[1]{\begin{lemma}{#1} \sl}
	\newcommand{\elem}{\end{lemma}}

\newenvironment{warning}[1]%
{ \begin{Warning} \label{W:#1}  \ifShowLabels \TeXref{W:#1} \fi }%
	{ \end{Warning} }
\newcommand{\war}[1]{\begin{warning}{#1} \sl}
	\newcommand{\ewar}{\end{warning}}

\newenvironment{assumption}[1]%
{ \begin{Assumption} \label{A:#1}  \ifShowLabels \TeXref{A:#1} \fi }%
	{ \end{Assumption} }
\newcommand{\ass}[1]{\begin{assumption}{#1} \sl}
	\newcommand{\eass}{\end{assumption}}	
	
\newenvironment{propos}[1]%
{ \begin{Prop} \label{P:#1}  \ifShowLabels \TeXref{P:#1} \fi }%
	{ \end{Prop} }
\newcommand{\prop}[1]{\begin{propos}{#1}\sl }
	\newcommand{\eprop}{\end{propos}}

\newenvironment{corol}[1]%
{ \begin{Cor} \label{C:#1}  \ifShowLabels \TeXref{C:#1} \fi }%
	{ \end{Cor} }
\newcommand{\cor}[1]{\begin{corol}{#1} \sl }
	\newcommand{\ecor}{\end{corol}}

\newenvironment{defeni}[1]%
{ \begin{Def} \label{D:#1}  \ifShowLabels \TeXref{D:#1} \fi }%
	{ \end{Def} }
\newcommand{\defe}[1]{\begin{defeni}{#1} \sl }
	\newcommand{\edefe}{\end{defeni}}

\newenvironment{remark}[1]%
{ \begin{Rem} \label{R:#1}  \ifShowLabels \TeXref{R:#1} \fi }%
	{ \end{Rem} }
\newcommand{\rem}[1]{\begin{remark}{#1}}
	\newcommand{\erem}{\end{remark}}

\newenvironment{conjec}[1]%
{ \begin{Conj} \label{Co:#1}  \ifShowLabels \TeXref{Co:#1} \fi }%
	{ \end{Conj} }
\renewcommand{\conj}[1]{\begin{conjec}{#1} \sl }
	\newcommand{\econj}{\end{conjec}}

\newenvironment{First proof}[1]%
{ \begin{First proof} \label{Co:#1}  \ifShowLabels \TeXref{Co:#1} \fi }%
	{ \end{First proof} }

\newenvironment{Second proof}[1]%
{ \begin{Second proof} \label{Co:#1}  \ifShowLabels \TeXref{Co:#1} \fi }%
	{ \end{Second proof} }

\newcommand{\eq}[1]%
{ \ifShowLabels \TeXref{E:#1} \fi
	\begin{equation} \label{E:#1} }
\newcommand{\eeq}{ \end{equation} }

\newcommand{\prf}{ \begin{proof} }
	\newcommand{\epr}{ \end{proof} }

\newcommand{\prft}{ \begin{proof} }
	\newcommand{\eprt}{ \end{proof} }

\newcommand\nc{\newcommand}

\nc{\HC}{{\mathcal{HC}}}
\nc{\on}{\operatorname}
\nc{\BA}{{\mathbb{A}}}
\nc{\BC}{{\mathbb{C}}}
\nc{\BF}{{\mathbb{F}}}
\nc{\BG}{{\mathbb{G}}}
\nc{\BM}{{\mathbb{M}}}
\nc{\BN}{{\mathbb{N}}}
\nc{\BO}{{\mathbb{O}}}
\nc{\BQ}{{\mathbb{Q}}}
\nc{\BP}{{\mathbb{P}}}
\nc{\BR}{{\mathbb{R}}}
\nc{\BZ}{{\mathbb{Z}}}
\nc{\BS}{{\mathbb{S}}}

\nc{\CA}{{\mathcal{A}}}
\nc{\CB}{{\mathcal{B}}}
\nc{\CalC}{{\mathcal C}}
\nc{\CalD}{{\mathcal D}}
\nc{\CE}{{\mathcal{E}}}
\nc{\CF}{{\mathcal{F}}}
\nc{\CG}{{\mathcal{G}}}
\nc{\CH}{{\mathcal{H}}}
\nc{\CK}{{\mathcal{K}}}
\nc{\CL}{{\mathcal{L}}}
\nc{\CM}{{\mathcal{M}}}
\nc{\CMM}{{\mathcal{M}^{\operatorname{gen}}_\hbar(-\rho)}}
\nc{\CN}{{\mathcal{N}}}
\nc{\CO}{{\mathcal{O}}}
\nc{\CP}{{\mathcal{P}}}
\nc{\CQ}{{\mathcal{Q}}}
\nc{\CR}{{\mathcal{R}}}
\nc{\CS}{{\mathcal{S}}}
\nc{\CT}{{\mathcal{T}}}
\nc{\CU}{{\mathcal{U}}}
\nc{\CV}{{\mathcal{V}}}
\nc{\CW}{{\mathcal{W}}}
\nc{\CX}{{\mathcal{X}}}
\nc{\CY}{{\mathcal{Y}}}
\nc{\CZ}{{\mathcal{Z}}}

\nc{\gen}{{\operatorname{gen}}}
\nc{\cM}{{\check{\mathcal M}}{}}
\nc{\csM}{{\check{\mathcal A}}{}}
\nc{\obM}{{\overset{\circ}{\mathbf M}}{}}
\nc{\oCA}{{\overset{\circ}{\mathcal A}}{}}
\nc{\obA}{{\overset{\circ}{\mathbf A}}{}}
\nc{\ooM}{{\overset{\circ}{M}}{}}
\nc{\osM}{{\overset{\circ}{\mathsf M}}{}}
\nc{\vM}{{\overset{\bullet}{\mathcal M}}{}}
\nc{\nM}{{\underset{\bullet}{\mathcal M}}{}}
\nc{\obD}{{\overset{\circ}{\mathbf D}}{}}
\nc{\cp}{{\overset{\circ}{\mathbf p}}{}}
\nc{\ofZ}{{\overset{\circ}{\mathfrak Z}}{}}

\nc{\fa}{{\mathfrak{a}}}
\nc{\fb}{{\mathfrak{b}}}
\nc{\fg}{{\mathfrak{g}}}
\nc{\fgl}{{\mathfrak{gl}}}
\nc{\fh}{{\mathfrak{h}}}
\nc{\fj}{{\mathfrak{j}}}
\nc{\fm}{{\mathfrak{m}}}
\nc{\fn}{{\mathfrak{n}}}
\nc{\fu}{{\mathfrak{u}}}
\nc{\fp}{{\mathfrak{p}}}
\nc{\frr}{{\mathfrak{r}}}
\nc{\fs}{{\mathfrak{s}}}
\nc{\ft}{{\mathfrak{t}}}
\nc{\fT}{{\mathfrak{T}}}
\nc{\ofT}{{\overline{\mathfrak T}}}
\nc{\ofS}{{\overline{\mathfrak S}}}
\nc{\fsl}{{\mathfrak{sl}}}
\nc{\hsl}{{\widehat{\mathfrak{sl}}}}
\nc{\hgl}{{\widehat{\mathfrak{gl}}}}
\nc{\hg}{{\widehat{\mathfrak{g}}}}
\nc{\chg}{{\widehat{\mathfrak{g}}}{}^\vee}
\nc{\hn}{{\widehat{\mathfrak{n}}}}
\nc{\chn}{{\widehat{\mathfrak{n}}}{}^\vee}

\nc{\fA}{{\mathfrak{A}}}
\nc{\fB}{{\mathfrak{B}}}
\nc{\fD}{{\mathfrak{D}}}
\nc{\fE}{{\mathfrak{E}}}
\nc{\fF}{{\mathfrak{F}}}
\nc{\fG}{{\mathfrak{G}}}
\nc{\fI}{{\mathfrak{I}}}
\nc{\fJ}{{\mathfrak{J}}}
\nc{\fK}{{\mathfrak{K}}}
\nc{\fL}{{\mathfrak{L}}}
\nc{\fM}{{\mathfrak{M}}}
\nc{\fN}{{\mathfrak{N}}}
\nc{\frP}{{\mathfrak{P}}}
\nc{\fQ}{{\mathfrak Q}}
\nc{\fS}{{\mathfrak S}}
\nc{\fU}{{\mathfrak{U}}}
\nc{\fZ}{{\mathfrak{Z}}}

\nc{\ba}{{\mathbf{a}}}
\nc{\bb}{{\mathbf{b}}}
\nc{\bc}{{\mathbf{c}}}
\nc{\bd}{{\mathbf{d}}}
\nc{\be}{{\mathbf{e}}}
\nc{\bi}{{\mathbf{i}}}
\nc{\bj}{{\mathbf{j}}}
\nc{\bn}{{\mathbf{n}}}
\nc{\bp}{{\mathbf{p}}}
\nc{\br}{{\mathbf{r}}}
\nc{\bq}{{\mathbf{q}}}
\nc{\bu}{{\mathbf{u}}}
\nc{\bv}{{\mathbf{v}}}
\nc{\bx}{{\mathbf{x}}}
\nc{\by}{{\mathbf{y}}}
\nc{\bw}{{\mathbf{w}}}
\nc{\bA}{{\mathbf{A}}}
\nc{\bB}{{\mathbf{B}}}
\nc{\bC}{{\mathbf{C}}}
\nc{\bD}{{\mathbf{D}}}
\nc{\bE}{{\mathbf{E}}}
\nc{\bK}{{\mathbf{K}}}
\nc{\bH}{{\mathbf{H}}}
\nc{\bM}{{\mathbf{M}}}
\nc{\bN}{{\mathbf{N}}}
\nc{\bO}{{\mathbf{O}}}
\nc{\bQ}{{\mathbf Q}}
\nc{\bS}{{\mathbf{S}}}
\nc{\bT}{{\mathbf{T}}}
\nc{\bV}{{\mathbf{V}}}
\nc{\bW}{{\mathbf{W}}}
\nc{\bX}{{\mathbf{X}}}
\nc{\bP}{{\mathbf{P}}}
\nc{\bZ}{{\mathbf{Z}}}

\nc{\sA}{{\mathsf{A}}}
\nc{\sB}{{\mathsf{B}}}
\nc{\sC}{{\mathsf{C}}}
\nc{\sD}{{\mathsf{D}}}
\nc{\sF}{{\mathsf{F}}}
\nc{\sK}{{\mathsf{K}}}
\nc{\sM}{{\mathsf{M}}}
\nc{\sO}{{\mathsf{O}}}
\nc{\sQ}{{\mathsf{Q}}}
\nc{\sP}{{\mathsf{P}}}
\nc{\sV}{{\mathsf{V}}}
\nc{\sW}{{\mathsf{W}}}
\nc{\sZ}{{\mathsf{Z}}}
\nc{\sfp}{{\mathsf{p}}}
\nc{\sr}{{\mathsf{r}}}
\nc{\st}{{\mathsf{t}}}
\nc{\sfb}{{\mathsf{b}}}
\nc{\sfc}{{\mathsf{c}}}
\nc{\sd}{{\mathsf{d}}}
\nc{\sg}{{\mathsf{g}}}
\nc{\sk}{{\mathsf{k}}}
\nc{\sfl}{{\mathsf{l}}}

\nc{\BK}{{\bar{K}}}

\nc{\tA}{{\widetilde{\mathbf{A}}}}
\nc{\tB}{{\widetilde{\mathcal{B}}}}
\nc{\tg}{{\widetilde{\mathfrak{g}}}}
\nc{\tG}{{\widetilde{G}}}
\nc{\TM}{{\widetilde{\mathbb{M}}}{}}
\nc{\tN}{{\widetilde{\mathcal{N}}}{}}
\nc{\tO}{{\widetilde{\mathsf{O}}}{}}
\nc{\tU}{{\widetilde{\mathfrak{U}}}{}}
\nc{\TZ}{{\tilde{Z}}}
\nc{\tZ}{\widetilde{Z}{}}
\nc{\tx}{{\tilde{x}}}
\nc{\tbv}{{\tilde{\bv}}}
\nc{\tfP}{{\widetilde{\mathfrak{P}}}{}}
\nc{\tz}{{\tilde{\zeta}}}
\nc{\tmu}{{\tilde{\mu}}}

\nc{\td}{\ddot{\underline{d}}{}}
\nc{\tzeta}{\widetilde{\zeta}{}}
\nc{\hd}{{\widehat{\underline{d}}}}
\nc{\hG}{{\widehat{G}}}
\nc{\hBP}{\widehat{\mathbb P}{}}
\nc{\hQ}{{\widehat{Q}}}
\nc{\hsM}{\widehat{\mathsf M}{}}
\nc{\hfM}{\widehat{\mathfrak M}{}}
\nc{\hCP}{\widehat{\mathcal P}{}}
\nc{\hCR}{\widehat{\mathcal R}{}}
\nc{\hCS}{{\widehat{\mathcal S}}}
\nc{\hfZ}{\widehat{\mathfrak Z}{}}
\nc{\hZ}{\widehat{Z}{}}

\nc{\urho}{\underline{\rho}}
\nc{\uB}{\underline{B}}
\nc{\uC}{{\underline{\mathbb{C}}}}
\nc{\ui}{\underline{i}}
\nc{\ofP}{{\overline{\mathfrak{P}}}}

\nc{\hrho}{{\hat{\rho}}}

\nc{\unl}{\underline}
\nc{\ol}{\overline}
\nc{\one}{{\mathbf{1}}}
\nc{\two}{{\mathbf{t}}}

\nc{\Sym}{{\mathop{\operatorname{Sym}}}}
\nc{\Tot}{{\mathop{\operatorname{\normalshape Tot}}}}
\nc{\Hilb}{{\mathop{\operatorname{\normalshape Hilb}}}}
\nc{\Hom}{{\mathop{\operatorname{Hom}}}}
\nc{\CHom}{{\mathop{\operatorname{{\mathcal{H}}\it om}}}}
\nc{\defi}{{\mathop{\operatorname{\normalshape def}}}}
\nc{\length}{{\mathop{\operatorname{\normalshape length}}}}

\nc{\Cliff}{{\mathsf{Cliff}}}
\nc{\Fl}{{\mathcal{F}\ell}}
\nc{\Fib}{{\mathsf{Fib}}}
\nc{\Coh}{{\mathsf{Coh}}}
\nc{\FCoh}{{\mathsf{FCoh}}}

\nc{\reg}{{\text{\normalshape reg}}}
\nc{\res}{{\operatorname{res}}}
\nc{\cplus}{{\mathbf{C}_+}}
\nc{\cminus}{{\mathbf{C}_-}}
\nc{\cthree}{{\mathbf{C}_*}}
\nc{\Qbar}{{\bar{Q}}}

\nc{\bh}{{\bar{h}}}
\nc{\bOmega}{{\overline{\Omega}}}
\nc\tGr{\widetilde{\Gr}}

\nc{\seq}[1]{\stackrel{#1}{\sim}}
\nc\ogu{\overline{G/U}}
\nc\chlam{\check{\lam}}

\nc\St{\operatorname{St}}
\nc{\oZ}{{\overset{\circ}{Z}}}
\nc{\tF}{\widetilde{\mathcal F}}

\nc\uS{\underline{S}}
\nc\QM{\mathcal{QM}}

\nc{\chmu}{\check{\mu}}
\newcommand\iso{\,\vphantom{j^{X^2}}\smash{\overset{\sim}{\vphantom{\rule{0pt}{0.20em}}\smash{\longrightarrow}}}\,}
\nc{\ul}{\underline}
\nc{\Mvd}{\mathfrak{M}(\underline{v},\underline{d})}
\nc{\MvdT}{\mathfrak{M}(\underline{v}^{\dagger},\underline{d}^{\dagger})}
\nc{\MVD}{\mathfrak{M}(V,D)}
\nc{\mt}{\mapsto}
\nc{\sm}{\setminus}
\nc{\ra}{\rightarrow}
\nc{\lar}{\leftarrow}
\nc{\hr}{\hookrightarrow}
\nc{\La}{\Lambda}
\nc{\Lap}{\Lambda^{+}}
\nc{\oZal}{\overset{\circ}{Z^{\alpha}}}
\nc{\sig}{\sigma}
\nc{\al}{\alpha}
\nc{\la}{\lambda}
\nc{\is}{\simeq}
\nc{\ip}{\iota^{+}_{\la, \mu}}
\nc{\im}{\iota^{-}_{\la, \mu}}
\nc{\jp}{j^{+}_{\la, \mu}}
\nc{\jm}{j^{-}_{\la, \mu}}
\nc{\pip}{\pi^{+}_{\la, \mu}}
\nc{\pim}{\pi^{-}_{\la, \mu}}
\nc{\s}{\star}
\nc{\fpt}{[A^{\la},B^{\la},\gamma^{\la},\delta^{\la}]}
\nc{\ulfpt}{[\ul{A}^{\la},\ul{B}^{\la},\ul{\gamma}^{\la},\ul{\delta}^{\la}]}
\nc{\lvee}{\!\scriptscriptstyle\vee}
\nc{\Gr}{{\operatorname{Gr}}}
\nc{\rra}{\twoheadrightarrow}

\nc{\End}{\on{End}}

\nc{\RHom}{R\mathcal{H}om}

\nc{\yg}{Y(\fg)}
\nc{\yvg}{Y_V(\fg)}
\nc{\CAg}{\CA_{\fg}}
\nc{\Ag}{A_{\fg}}
\nc{\Sgt}{S^{\bullet}(\fg[t])}
\nc{\Sg}{S^{\bullet}(\fg)}
\nc{\Ugt}{U(\fg[t])}

\nc{\Spec}{\operatorname{Spec}}

\author{Vasily Krylov}
\address{V.K.: Department of Mathematics
Massachusetts Institute of Technology
\newline
77 Massachusetts Avenue,
Cambridge, MA 02139,
United States of America;
\newline National Research University Higher School of Economics, Russian Federation\newline
Department of Mathematics, 6 Usacheva st., Moscow 119048;
}
\email{krvas@mit.edu, krylovasya@gmail.com}
\author{Leonid Rybnikov}
\address{L.R.: Department of Mathematics
Massachusetts Institute of Technology
\newline
77 Massachusetts Avenue,
Cambridge, MA 02139,
United States of America;
\newline On leave from National Research University Higher School of Economics, Russian Federation\newline
Department of Mathematics, 6 Usacheva st., Moscow 119048;
}
\email{lrybnikov@math.harvard.edu, leo.rybnikov@gmail.com}

\begin{document}
	\begin{abstract}
The loop group $G((z^{-1}))$ of a simple complex Lie group $G$ has a natural Poisson structure. We introduce a natural family of Poisson commutative subalgebras $\ol{{\bf{B}}}(C) \subset \CO(G((z^{-1}))$ depending on the parameter $C\in G$ called {\emph{classical universal Bethe subalgebras}}.
To every antidominant cocharacter $\mu$ of the maximal torus $T \subset G$ one can associate the closed Poisson subspace $\CW_\mu$ of $G((z^{-1}))$ (the Poisson algebra $\CO(\CW_\mu)$ is the classical limit of so-called shifted Yangian $Y_\mu(\mathfrak{g})$ defined in \cite[Appendix B]{bfn}). We consider the images of $\ol{{\bf{B}}}(C)$ in $\CO(\CW_\mu)$, denoted by $\ol{B}_\mu(C)$, that should be considered as classical versions of (not yet defined in general) Bethe subalgebras in shifted Yangians. For regular $C$ centralizing $\mu$, we compute the Poincar\'e series of these subalgebras.  
For $\mathfrak{g}=\mathfrak{gl}_n$, we define the natural quantization ${\bf{Y}}^{\mathrm{rtt}}(\mathfrak{gl}_n)$ of $\CO(\on{Mat}_n((z^{-1})))$
	and universal Bethe subalgebras ${\bf{B}}(C) \subset {\bf{Y}}^{\mathrm{rtt}}(\mathfrak{gl}_n)$. Using the RTT realization of $Y_\mu(\mathfrak{gl}_n)$ (invented by Frassek, Pestun, and Tsymbaliuk), we obtain the natural surjections ${\bf{Y}}^{\mathrm{rtt}}(\mathfrak{gl}_n) \twoheadrightarrow Y_\mu(\mathfrak{gl}_n)$ which quantize the embedding $\CW_\mu \subset \on{Mat}_n((z^{-1}))$. Taking the images of ${\bf{B}}(C)$ in $Y_\mu(\mathfrak{gl}_n)$ we recover Bethe subalgebras $B_\mu(C) \subset Y_\mu(\mathfrak{gl}_n)$ proposed by Frassek, Pestun and Tsymbaliuk. 
		\end{abstract}

	\title[]{Bethe subalgebras in antidominantly shifted Yangians}
	
	\maketitle

\section{Introduction}


\subsection{Yangians, Bethe subalgebras and their classical limits}{}  Let $\mathfrak{g}$ be a simple finite dimensional Lie algebra over complex numbers and $G$ be the corresponding simply connected algebraic group. The Yangian $Y(\mathfrak{g})$ is the Hopf algebra deformation of the algebra $\CO(G[[z^{-1}]]_1)$ of functions on the first congruence subgroup deforming the natural Poisson structure on the group scheme $G[[z^{-1}]]_1$ (see \cite{kwwy}). In other words, there exists a filtration on $Y(\mathfrak{g})$ such that $\on{gr}Y(\mathfrak{g}) \simeq \CO(G[[z^{-1}]]_1)$ as Poisson algebras (see for example \cite[Section 2]{ir}). 

Bethe subalgebras are a family of commutative subalgebras $B(C) \subset Y(\mathfrak{g})$ depending on $C \in G$. In \cite[Section 4]{ir} the associated graded algebra $\ol{B}(C)=\on{gr}B(C) \subset \CO(G[[z^{-1}]]_1)$ was described.  

In \cite{bk}, \cite{kwwy}, \cite[Appendix B]{bfn} generalizations of $Y(\mathfrak{g})$ were defined. Let us fix a Borel subgroup $B \subset G$ and a maximal torus $T \subset B$. Let $\La$ be the cocharacter lattice of $T$.  To every $\mu \in \La$, one can associate the so-called shifted Yangian $Y_\mu(\mathfrak{g})$. By the results of  \cite[Section 5]{fkprw}, the algebra $Y_\mu(\mathfrak{g})$ is filtered and the associated graded Poisson algebra $\on{gr}Y_\mu(\mathfrak{g})$ can be identified with the algebra of functions on a Poisson scheme that the authors of \cite{fkprw} denote by $\CW_\mu$. This scheme is a closed subscheme of $G((z^{-1}))$. It follows from \cite[Theorem A.8]{kpw} that for antidominant $\mu$ the embedding $\CW_\mu \subset G((z^{-1}))$ is {\emph{Poisson}} so that the Poisson algebra $\CO(\CW_\mu)$ is the quotient of the Poisson algebra $\CO(G((z^{-1})))$ by the Poisson ideal.

\rem{}
{\emph{The antidominant $\mu$ is distinguished since for such $\mu$ we have 
\begin{equation*}
\CW_\mu=G[[z^{-1}]]_1 z^{\mu} G[[z^{-1}]]_1 \subset G((z^{-1})),
\end{equation*}
where $z^\mu \in T((z^{-1})) \subset G((z^{-1}))$ is the element corresponding to $\mu\colon \on{Spec}\BC[z^{\pm 1}] \ra T$.}}
\erem

Recall now that for $\mu=0$ we have a family of commutative subalgebras $B(C) \subset Y_0(\mathfrak{g})=Y(\mathfrak{g})$ (depending on $C \in G$) and their ``classical'' limits $\ol{B}(C) \subset \CO(\CW_0)$. The natural question is the following. 

\quad

{\bf{Question}}: can we define a family of  commutative (resp. Poisson commutative) subalgebras 
\begin{equation*}
B_\mu(C) \subset Y_\mu(\mathfrak{g}),\, \ol{B}_\mu(C) \subset \CO(\CW_\mu)    
\end{equation*}
for antidominant $\mu \in \La$ generalizing families $B(C)$, $\ol{B}(C)$ above?

\quad

In this note, we answer this question for algebras $\ol{B}_\mu(C)$. For $\mathfrak{g}=\mathfrak{gl}_n$ algebras $B_\mu(C)$ were defined by Frassek, Pestun and Tsymbaliuk (private communication) using 
their RTT realization of antidominantly shifted Yangians (see \cite{fpt}). Let us recall the definition of $B_\mu(C)$ (see Section \ref{bethe_subalg_in_shifted_section!!!} for details). The Yangian $Y_\mu^{\mathrm{rtt}}(\mathfrak{gl}_n)$ is generated by $\{t_{ij}^{(r)}\}^{r \in \BZ}_{1 \leqslant i,j \leqslant n}$ (see \cite[Section 2.3]{fpt} or Definition \ref{defe_RTT_Yang!} for details). For $i,j = 1,\ldots,n$ we set $t_{ij}(u):=\sum_{r \in \BZ}t_{ij}^{(r)}u^{-r} \in Y_\mu^{\mathrm{rtt}}(\mathfrak{gl}_n)((u^{-1}))$.
Let $T(u) \in Y_\mu^{\mathrm{rtt}}(\mathfrak{gl}_n)((u^{-1})) \otimes \on{End}(\BC^n)$ be the matrix $T(u):=(t_{ij}(u))_{ij}$. 
For $k=1,\ldots,n$ we denote by $A_k \in \on{End}(\BC^n)^{\otimes k}$ the antisymmetrization  map normalized so that $A_k^2=A_k$. We set
\begin{equation*}
\tau_{\mu,k}(u,C):=\on{tr}_{(\BC^n)^{\otimes k}} A_k C_1 \ldots C_k T_1(u) \ldots T_k(u-k+1) \in Y^{\mathrm{rtt}}_\mu(\mathfrak{gl}_n)((u^{-1})).
\end{equation*}
The algebra $B_\mu(C)$ is generated by the coefficients of $\tau_{\mu,k}(u,C),\,k=1,\ldots,n$.
We denote the coefficient of $u^{-r}$ in $\tau_{\mu,k}(u,C)$ by $\tau_{\mu,k}(C)^{(r)}$. Set $\omega_k^*:=-\epsilon_n^\vee-\ldots-\epsilon_{n-k+1}^\vee \in (\BC^n)^*$, where $\epsilon_1^{\vee},\ldots,\epsilon_n^\vee$ is the standard basis of $(\BC^n)^*$. The authors of \cite{fpt}  formulated the following conjecture  about the structure of the algebras $B_\mu(C)$ for regular $C$:

\begin{conjec}{}[Frassek, Pestun, and Tsymbaliuk, 2019]\label{FPT_conj_form}

$A)$ For $C \in \on{GL}_n^{\mathrm{reg}}$ the algebra $B_\mu(C)$ is a polynomial algebra in the elements $\{ \tau^{(r)}_{\mu,k}(C) \,|\, r > \langle \omega_k^*,\mu \rangle\}$.

$B)$ For $C \in  \on{GL}_n^{\mathrm{reg}}$ the subalgebra $B_\mu(C) \subset Y_\mu(\mathfrak{g})$ is maximal commutative.
\end{conjec}

\rem{}
\em{Conjecture \ref{FPT_conj_form} for $\mu=0$ and $C \in T^{\mathrm{reg}}$ follows from \cite[Theorem 1.3]{no}.}
\erem

In this note, we prove part $A)$ of the Conjecture \ref{FPT_conj_form}  
(assuming that $C$ lies in the centralizer of $\mu$ in $G$). We prove that for $C \in Z_G(\mu)^{\mathrm{reg}}$ the algebra $B_\mu(C)$ is a polynomial algebra in the elements $\{ \tau^{(r)}_{\mu,k}(C) \,|\, r > \langle \omega_k^*,\mu \rangle\}$. In particular, it follows that for $\mu$ such that $Z_G(\mu)=T$ claim of Conjecture \ref{FPT_conj_form} A) holds for {\em{any}} $C \in T$.


\subsection{Main results of the paper}{} 
\subsubsection{Classical situation} 
Let us start with the classical situation. Recall that $\mu \in \La$ is antidominant. Recall that  $\CO(\CW_\mu)$ is the quotient of the bigger Poisson algebra $\CO(G((z^{-1})))$. It turns out (see Definition \ref{defe_univ_class_bethe} and Proposition \ref{comm_univ_class_bethe}) that one can define a family $\ol{{\bf{B}}}(C) \subset \CO(G((z^{-1})))$, $C \in G$, of Poisson commutative subalgebras (that we call classical universal Bethe subalgebras) and then  $\ol{B}_\mu(C) \subset \CO(\CW_\mu)$ can be defined as the image of $\ol{\bf{B}}(C)$. The algebra $\ol{{\bf{B}}}(C)$ is the subalgebra of $\CO(G((z^{-1})))$ generated by the Fourier coefficients of the functions 
\begin{equation*}
G((z^{-1})) \ni g \mapsto \on{tr}_{V_{\omega_i}} \rho_i(C)\rho_i(g) \in \BC((z^{-1})),
\end{equation*}
where $\rho_i\colon G \ra \on{End}(V_{\omega_i})$ are the fundamental representations of $G$.

The scheme $\CW_\mu$ has a natural $\BC^\times$-action that induces the grading on $\CO(\CW_\mu)$. The subalgebra $\ol{B}_\mu(C) \subset \CO(\CW_\mu)$ is not graded in general, but it has the induced filtration. Let $L:=Z_{G}(\mu)$ be the centralizer of $\mu$ in $G$. For $C \in L$ we denote by $\ol{B}_L(C) \subset \CO(L[[z^{-1}]]_1)$ the (classical) Bethe subalgebra of $\CO(L[[z^{-1}]]_1)$. 
Consider the closed embedding $L[[z^{-1}]]_1 \subset \CW_\mu$ given by $g \mapsto gz^\mu$. This embedding induces  surjection $\CO(\CW_\mu) \twoheadrightarrow \CO(L[[z^{-1}]]_1)$ at the level of functions.
The main result of Section \ref{class-cal_bethe_section} is the following theorem (see Theorem \ref{size_reg_bethe} for  details).

\begin{thmx}\label{thm_A}
If $C \in L^{\mathrm{reg}}$ then the composition 
$
\on{gr}\ol{B}_\mu(C) \hookrightarrow \CO(\CW_\mu) \twoheadrightarrow \CO(L[[z^{-1}]]_1)
$
induces an isomorphism $\on{gr}\ol{B}_\mu(C) \iso \ol{B}_L(C)$.
\end{thmx}
As a corollary, we conclude that the algebra $\ol{B}_\mu(C)$ is a free polynomial algebra, and the size of $\ol{B}_\mu(C)$ is the same as the size of $\ol{B}_{L}(C)$ (see Corollary \ref{cor_size_class_Bethe_shifted}).

\rem{} 
{\emph{By ``size'' we always  mean a Poincar\'e series with respect to a certain natural filtration.}}
\erem

Let us briefly describe the idea of the proof of Theorem \ref{thm_A}. The natural generators of $\ol{B}_\mu(C)$ are the Fourier coefficients of the functions $g \mapsto \on{tr}_{V_{\omega_i}} \rho_i(C)\rho_i(g)$. Taking their images in $\CO(L[[z^{-1}]]_1)$, one can see that the highest component of this trace is nothing else but the trace of the same operator but restricted to the subspace of $V_{\omega_i}$ generated by vectors of $V_{\omega_i}$ having weight $\nu$ such that $\langle \nu,\mu \rangle = \langle w_0(\omega_i),\mu \rangle$ (here $w_0$ is the longest element of the Weyl group of $G$). This vector subspace is nothing else but the irreducible representation of $L$ with the lowest weight $w_0(\omega_i)$. Now Theorem \ref{thm_A} can be deduced from \cite[Proposition 4.6]{ir}.

\subsubsection{Quantum situation} Let us now consider the ``quantum'' situation. In other words, our goal is to construct Bethe subalgebras in the shifted Yangians $Y_\mu(\mathfrak{g})$ (see \cite[Definition B.2]{bfn}, or \cite[Definition 3.5]{fkprw}). The only known definition of Bethe subalgebras in the standard Yangian $Y(\mathfrak{g})$ uses the so-called RTT realization of $Y(\mathfrak{g})$. Such a realization of $Y_{\mu}(\mathfrak{g})$ is not known in general but was obtained in \cite{fpt}, \cite{ft} for $\mathfrak{g}$ of type $ABCD$. In this note, we restrict our discussion only to the type $A$ case. 
Using the RTT realization of $Y_\mu(\mathfrak{gl}_n)$ we recall the definition of Bethe subalgebras $B_\mu(C) \subset Y_\mu(\mathfrak{gl}_n)$ which belongs to \cite{fpt} (see Definition \ref{defe_bethe_shifted}). 
We actually show that the algebras $B_\mu(C)$ can be considered as the images of one ``universal'' Bethe subalgebra ${\bf{B}}(C)$ inside some algebra ${\bf{Y}}^{\mathrm{rtt}}(\mathfrak{gl}_n)$ quantizing $\CO(\on{Mat}_n((z^{-1})))$  similarly to how all $\overline{B}_\mu(C)$ are images of a single $\overline{\bf{B}}(C)\subset \CO(\on{Mat}_n((z^{-1})))$ (see Section \ref{univ_bethe_yang_sect!} for details).

\begin{thmx}\label{thm_B} The
Poisson algebra $\CO(\on{Mat}_n((z^{-1})))$ has a natural quantization  ${\bf{Y}}^{\mathrm{rtt}}(\mathfrak{gl}_n)$ that maps surjectively onto $Y_\mu(\mathfrak{gl}_n)$ for every antidominant $\mu \in \La$. The subalgebra $B_\mu(C) \subset Y_\mu(\mathfrak{gl}_n)$ is the image of a ``universal'' Bethe subalgebra ${\bf{B}}(C) \subset {\bf{Y}}^{\mathrm{rtt}}(\mathfrak{gl}_n)$.
\end{thmx}

\rem{}
{\em{The surjection ${\bf{Y}}^{\mathrm{rtt}}(\mathfrak{gl}_n) \twoheadrightarrow Y_\mu(\mathfrak{gl}_n)$ is a quantization of the restriction homomorphism $\CO(\on{Mat}_n((z^{-1}))) \twoheadrightarrow \CO(\CW_\mu)$, see Corollary \ref{surj_quant_surj} for details.}}
\erem


Let us briefly describe the construction of the algebras ${\bf{Y}}^{\mathrm{rtt}}(\mathfrak{gl}_n)$, ${\bf{B}}(C)$. Recall that the standard RTT Yangian $Y(\mathfrak{gl}_n)$ has the following realization. As an algebra over $\BC$ it is generated by $\{t_{ij}^{(r)}\}^{r \in \BZ_{\geqslant 0}}_{1 \leqslant i,j \leqslant n}$ subject to the following relations:
\begin{equation}\label{eq_def_1_stand}
[t_{ij}^{(p+1)},t_{kl}^{(q)}]-[t_{ij}^{(p)},t_{kl}^{(q+1)}]=t_{kj}^{(q)}t_{il}^{(p)}-t_{kj}^{(p)}t_{il}^{(q)},\, p,q \in \BZ_{\geqslant 0},  
\end{equation}
\begin{equation*}
t_{ij}^{(0)}=\delta_{ij}.    
\end{equation*}
We define ${\bf{Y}}^{\mathrm{rtt}}(\mathfrak{gl}_n)^{\mathrm{pol}}$ as the algebra over $\BC$ generated by $\{t_{ij}^{(r)}\}^{r \in \BZ}_{1 \leqslant i,j \leqslant n}$ subject to the relations (\ref{eq_def_1_stand}) with $p,q \in \BZ$. For $N \in \BZ_{\geqslant 0}$ consider the two-sided ideal $I_N \subset {\bf{Y}}^{\mathrm{rtt}}(\mathfrak{gl}_n)^{\mathrm{pol}}$ generated by the elements 
$
\{t_{ij}^{(-r)}\,|\, r > N\}
$
and set 
${\bf{Y}}^{\mathrm{rtt}}(\mathfrak{gl}_n)_{N}:={\bf{Y}}^{\mathrm{rtt}}(\mathfrak{gl}_n)^{\mathrm{pol}}/I_N$, ${\bf{Y}}^{\mathrm{rtt}}(\mathfrak{gl}_n):=\underset{\longleftarrow}{\on{lim}}\,{\bf{Y}}^{\mathrm{rtt}}(\mathfrak{gl}_n)_N.$
The subalgebra ${\bf{B}}(C) \subset {\bf{Y}}^{\mathrm{rtt}}(\mathfrak{gl}_n)$ is defined using the same formulas as one uses in the definition of  standard Bethe subalgebra $B(C) \subset Y(\mathfrak{gl}_n)$ (see Section \ref{univ_bethe_in_yang} for details).

The main result of Section \ref{size_shift_bethe}  is the following theorem (see Corollary \ref{cor_size_Bethe_type_A} for details).

\begin{thmx}\label{thm_C}
For $C \in L^{\mathrm{reg}}$, the 
composition $\on{gr}B_\mu(C) \hookrightarrow \CO(\CW_\mu) \twoheadrightarrow \CO(L[[z^{-1}]]_1)$ induces the isomorphism $\on{gr}B_\mu(C) \iso \ol{B}_L(C)$. 
\end{thmx}

As a corollary, we conclude that the algebra  $B_\mu(C)$ is a free polynomial algebra and has the same size as the Cartan subalgebra $H \subset Y_\mu(\mathfrak{gl}_n)$ (see Corollary \ref{cor_size_Bethe_type_A}), which proves Conjecture \ref{FPT_conj_form} $A)$ for $C \in L^{\mathrm{reg}}$.


The idea of the proof of Theorem \ref{thm_C} is the following: we reduce it to Theorem \ref{thm_A} by showing that the associated graded of the natural generators of $B_\mu(C)$ and of $\ol{B}_\mu(C)$ are equal (considered as the elements of $\CO(\CW_\mu)$).

Overall, this note should be considered as an attempt to draw attention to study families of commutative subalgebras in  shifted Yangians and of their classical and ``universal'' versions.

\subsection{Relation with generalized transversal slices and truncated shifted Yangians}
Recall that $G$ is a simply connected Lie group corresponding to the simple Lie algebra $\mathfrak{g}$. Recall also that we fixed a maximal torus and a Borel subgroup $T \subset B \subset G$. To every dominant cocharacter $\la\colon \BC^\times \ra T$ and arbitrary $\mu\colon \BC^\times \ra T$ one can associate the closed Poisson subscheme $\ol{\CW}^\la_\mu \subset \CW_\mu$ defined as follows:
\begin{equation*}
\ol{\CW}^\la_\mu:= \CW_\mu \cap \ol{G[z]z^\la G[z]},    
\end{equation*}
where $G[z]$ is the space of maps $\on{Spec}(\BC[z]) \ra G$ and the closure is taken in $G((z^{-1}))$.
The scheme $\ol{\CW}^\la_\mu$ was introduced in the paper \cite{bfn} and is called a generalized transversal slice in the affine Grassmannian of $G$. 
Assume now that $\mu$ is antidominant.
Since the embedding $\ol{\CW}^\la_\mu \subset \CW_\mu$ is Poisson, we then can define a family of Poisson commutative subalgebras $\ol{B}^\la_\mu(C) \subset \CO(\ol{\CW}^\la_\mu)$ depending on $C \in G$. Indeed, we just need to take the images of $\ol{B}_\mu(C) \subset \CO(\CW_\mu)$ in $\CO(\ol{\CW}^\la_\mu)$.
The following question was raised by the authors of \cite{fpt}:
is it true that for $C \in T^{\mathrm{reg}}$ the subalgebra $\ol{B}^\la_\mu(C) \subset \CO(\ol{\CW}^\la_\mu)$ has transcendence degree equal to $\frac{1}{2}\on{dim}\ol{\CW}^{\la}_\mu$ (i.e. that the subalgebra $\ol{B}^{\la}_\mu(C) \subset \CO(\ol{\CW}^\la_\mu)$ defines a ``complete set of integrals'' in $\CO(\ol{\CW}^\la_\mu)$)? We plan to address this question (for generic $C$) in the future (the idea is to ``degenerate'' $\ol{B}_\mu(C)$ to the  subalgebra $\CO(T[[z^{-1}]]_1z^\mu) \subset \CO(\CW_\mu)$ whose image in  $\CO(\ol{\CW}^\la_\mu)$  has transcendence degree equal to $\frac{1}{2}\on{dim}\ol{\CW}^{\la}_\mu$).

The Poisson algebras $\CO(\ol{\CW}^\la_\mu)$ have natural quantizations $Y^\la_\mu(\mathfrak{g})$ called truncated shifted  Yangians (they were defined in \cite[Appendix B]{bfn} and in \cite[Theorem 4.7]{kpw} the authors proved that $Y^\la_\mu(\mathfrak{g})$ indeed quantizes $\CO(\ol{\CW}^\la_\mu)$). The algebra $Y^\la_\mu(\mathfrak{g})$ is the quotient of the shifted Yangian $Y_\mu(\mathfrak{g})$. Let us now restrict ourselves to the case $\mathfrak{g}=\mathfrak{gl}_n$. Then, taking images of $B_\mu(C) \subset Y_\mu(\mathfrak{gl}_n)$ in $Y^\la_\mu(\mathfrak{gl}_n)$, we obtain the family of commutative subalgebras $B^\la_\mu(C) \subset Y^\la_\mu(\mathfrak{gl}_n)$ depending on $C \in \on{GL}_n$.  It would be very interesting to study this family. Let us point out that the existence of this family as well as its interpretation through commutative subalgebras of quantized Coulomb branch was already noticed in the introduction of the paper \cite{fpt}. 

\subsection{Structure of the paper}{} In Section \ref{class-cal_bethe_section}, we recall the explicit description of the Poisson bracket on $\CO(G((z^{-1})))$, define the classical universal Bethe subalgebra $\ol{{\bf{B}}}(C) \subset \CO(G((z^{-1})))$ (see Definition \ref{defe_univ_class_bethe}), prove its Poisson commutativity (see Proposition \ref{comm_univ_class_bethe}) and then show that for $C \in G^{\mathrm{reg}}$ the natural generators of $\ol{{\bf{B}}}(C)$ are algebraically independent (see Proposition \ref{size_of_univ_class_Bethe}). We then recall closed Poisson subschemes $\CW_\mu \subset G((z^{-1}))$ and consider the classical Bethe subalgebras $\ol{B}_\mu(C) \subset \CO(\CW_\mu)$ (see Definition \ref{class_bethe_in_shifted}). We then prove Theorem \ref{thm_A}: for $C \in L^{\mathrm{reg}}$ (for example, if $C$ belongs to the regular part of a maximal torus $T \subset G$) we compute the size (Poincar\'e series) of the algebra $\ol{B}_\mu(C)$ and show that $\on{gr}\ol{B}_\mu(C)$ identifies naturally with the Bethe subalgebra $\ol{B}_{L}(C) \subset \CO(L[[z^{-1}]]_1)$ (see Theorem \ref{size_reg_bethe} and Corollary \ref{cor_size_class_Bethe_shifted}). 

In Section \ref{univ_bethe_yang_sect!}, we consider the case $\mathfrak{g}=\mathfrak{gl}_n$ and define the ``Yangian'' quantization of the Poisson algebra $\CO(\on{Mat}_n((z^{-1})))$ (see Definition \ref{defe_RTT_Yang_univ!}). We prove the analog of PBW theorem for this quantization ${\bf{Y}}^{\mathrm{rtt}}(\mathfrak{gl}_n)$ (see Proposition \ref{gen_quot_yang_univ}) and then show that $\on{gr}{\bf{Y}}^{\mathrm{rtt}}(\mathfrak{gl}_n) \simeq \CO(\on{Mat}_n((z^{-1})))$ as Poisson algebras (see Proposition \ref{iso_poiss_gr_yang_gl_univ}). We then define universal Bethe subalgebras ${\bf{B}}(C) \subset {\bf{Y}}^{\mathrm{rtt}}(\mathfrak{gl}_n)$ (see Definition \ref{defe_bethe_univ}). 
For $C \in G^{\mathrm{reg}}$ we show that $\on{gr}{\bf{B}}(C)=\ol{{\bf{B}}}(C)$ and conclude that  ${\bf{B}}(C)$ is a polynomial algebra in the natural generators (see Proposition \ref{gen_of_fat_B_prop}).

In Section \ref{shifted_yangians_section}, following \cite{fpt} we recall two definitions of  antidominantly shifted Yangians for $\mathfrak{gl}_n$ (see Definitions \ref{BFN_def_shifted_yangian}, \ref{defe_RTT_Yang!}) and then describe functions $\on{gr}t_{ij}^{(r)} \in \on{gr}Y^{\mathrm{rtt}}_\mu(\mathfrak{gl}_n)=\CO(\CW_\mu)$ (see Lemma \ref{gr_t_i_j!}). In Corollary \ref{surj_quant_surj}, 
we prove that the natural surjection ${\bf{Y}}^{\mathrm{rtt}}(\mathfrak{gl}_n) \twoheadrightarrow Y^{\mathrm{rtt}}_\mu(\mathfrak{gl}_n)$ quantizes the natural surjection $\CO(\on{Mat}_n((z^{-1}))) \twoheadrightarrow \CO(\CW_\mu)$.

In Section \ref{bethe_subalg_in_shifted_section!!!}, we recall Bethe subalgebras $B_\mu(C) \subset Y^{\mathrm{rtt}}_\mu(\mathfrak{gl}_n)$ (see Definition \ref{defe_bethe_shifted}) and remind their commutativity from the results of Section \ref{univ_bethe_in_yang} (see Proposition \ref{shifted_bethe_are_comm}). 
Theorem \ref{thm_B} then follows from the results of Sections \ref{univ_bethe_yang_sect!}, \ref{shifted_yangians_section}, \ref{bethe_subalg_in_shifted_section!!!}. We finish this note with the proof of Theorem \ref{thm_C}: assuming that $C \in L^{\mathrm{reg}}$ we describe the natural algebraically independent generators of $B_\mu(C)$, compute its Poincar\'e series and show that $\on{gr}B_\mu(C)$ is isomorphic to $\ol{B}_L(C)$ (see Theorem \ref{thm_size_shifted} and Corollary \ref{cor_size_Bethe_type_A}).

The note also contains two appendices. In Appendix \ref{app_1}, we prove that already for $\mathfrak{g}=\mathfrak{sl}_2$ (classical) Bethe subalgebras $\ol{B}_\mu(C) \subset \CO(\CW_\mu)$ can not be obtained as pullbacks of classical Bethe subalgebras in $\CO(\CW_0)$ (see Proposition \ref{alg_are_diff_prf}). Appendix \ref{app} contains a generalization of a well-known theorem of Steinberg that we are using in our arguments (see Proposition \ref{prop_lin_ind_inv}).




\subsection{Acknowledgements}{} We would like to thank  Michael Finkelberg, Aleksei Ilin, and Alexander Tsymbaliuk  for helpful
discussions and explanations. We are grateful to Alexander Tsymbaliuk for his useful comments on the preliminary version of the text.
We are also extremely grateful to our anonymous referees for the very careful proofreading of the text, for the suggestions on how to improve the exposition, and for pointing out many inaccuracies,  numerous typos, and historical errors.
Both of the authors were partially supported by the Foundation for the Advancement of Theoretical Physics and Mathematics ``BASIS''. L.R. is grateful to the Institut des Hautes Études Scientifiques and especially to Maxim Kontsevich for the hospitality and for the opportunity to avoid further political persecution in Russia and to continue working on this project.

\section{Classical ``universal'' Bethe subalgebra in $\CO(G((z^{-1})))$ and its images in $\CO(\CW_\mu)$}\label{class-cal_bethe_section}

\subsection{Classical universal Bethe subalgebra}{} Let $\mathfrak{g}$ be a finite dimensional  simple Lie algebra over $\BC$, we denote by  $G$ the simply connected group with Lie algebra $\mathfrak{g}$. 
Recall that we fix a Borel subgroup and a maximal torus 
$T \subset B \subset G$. Recall that  $\La=\on{Hom}(\BC^\times,T)$ is the cocharacter lattice of $T$. We denote by $I$ the set parametrizing vertices of the Dynkin diagram of $\mathfrak{g}$. For $i \in I$ we denote by $\omega_i$ the corresponding fundamental weight of $\mathfrak{g}$ and by $V_{\omega_i}$ the corresponding (fundamental) representation $\rho_i\colon \mathfrak{g} \ra \on{End}(V_{\omega_i})$, abusing the notation we denote by the same letter $\rho_i$ the corresponding representation of $G$.

Let $(\,,\,)$ be an invariant nondegenerate form on $\mathfrak{g}$ and let $\{x^a\}_{a=1,\ldots,\on{dim}\mathfrak{g}}$,  $\{x_a\}_{a=1,\ldots,\on{dim}\mathfrak{g}}$ be a pair of dual bases (of $\mathfrak{g}$ w.r.t. $(\,,\,)$). Let $V$ be a finite dimensional representation of $G$ and pick $v \in V, \beta \in V^*$. The matrix entry $\Delta_{\beta,v}(g)$ is the function on $G$ given by $\Delta_{\beta,v}(g)=\langle \beta,gv \rangle$.  Using this matrix entry we can define the function $\Delta^{(r)}_{\beta,v}$, $r \in \BZ$, on $G((z^{-1}))$ whose  value at $g(z^{-1}) \in G((z^{-1}))$ is the coefficient of $z^{-r}$ in $\Delta_{\beta,v}(g(z^{-1}))$. More precisely, these are given by the formula:
\begin{equation*}
\langle \beta, g(z^{-1})v \rangle =\sum_{r \in \BZ} \Delta^{(r)}_{\beta,v}(g(z^{-1})) z^{-r}.
\end{equation*}
It is convenient to introduce the following series ($u$ is a formal parameter): 
\begin{equation*}
\Delta_{\beta,v}(u):=\sum_{r \in \BZ} \Delta_{\beta,v}^{(r)} u^{-r} \in \CO(G((z^{-1})))[[u,u^{-1}]].
\end{equation*}

There is a nondegenerate pairing on $\mathfrak{g}((z^{-1}))$ coming from the residue and the invariant form $(\,,\,)$ on $\mathfrak{g}$. In particular,  $(\mathfrak{g}((z^{-1})),\mathfrak{g}[z],z^{-1}\mathfrak{g}[[z^{-1}]])$ is a Manin triple. This induces a Poisson-Lie group structure on $G((z^{-1}))$.

\rem{}
{\em{Let us recall that $G((z^{-1}))$
is and ind-scheme of ind-infinite type (see, for example, \cite[Proposition 2.5.1]{kv}). One way to see this is to fix a closed embedding $G \subset \BA^d$ for some $d \in \BZ_{\geqslant 0}$. This gives a closed embedding of $G((z^{-1}))$ inside the ind-scheme $\underset{\longrightarrow_N}{\on{lim}}\, \on{Spec}[a^{(l)}_i\,|\, i =1,\ldots,d,\, l \in \BZ_{\geqslant -N}]$. Preimages of $\on{Spec}[a^{(l)}_i\,|\, i =1,\ldots,d,\, l \in \BZ_{\geqslant -N}]$ realize $G((z^{-1}))$ as an ind-scheme.
}}
\erem

The Poisson bracket on $\CO(G((z^{-1})))$ is given by the following formula (see \cite[Proposition 2.13]{kwwy} and \cite[Remark A.7]{kpw}). Let $\rho_{V_1}\colon \mathfrak{g} \ra \on{End}(V_1)$, $\rho_{V_2}\colon \mathfrak{g} \ra \on{End}(V_2)$ be finite dimensional representations of $\mathfrak{g}$ (that can be also considered as representations of $G$). Pick $v_1 \in V_1$, $\beta_1 \in V_1^*$, $v_2 \in V_2$, $\beta_2 \in V_2^*$, $r,s \in \BZ$. Then we have 
\begin{equation*}
\{\Delta^{(r+1)}_{\beta_1,v_1},\Delta^{(s)}_{\beta_2,v_2}\}-\{\Delta^{(r)}_{\beta_1,v_1},\Delta^{(s+1)}_{\beta_2,v_2}\}=\sum_{a=1}^{\on{dim}\mathfrak{g}}  \Delta^{(r)}_{\beta_1,x_av_1}\Delta^{(s)}_{\beta_2,x^av_2} -\Delta^{(r)}_{x_a\beta_1,v_1}\Delta^{(s)}_{x^a\beta_2,v_2}      
\end{equation*}
that is equivalent to the following equality 
\begin{multline}\label{equal_poiss_in_vectors}
(u_1-u_2)\{\Delta_{\beta_1,v_1}(u_1),\Delta_{\beta_2,v_2}(u_2)\}=\\
=\sum_{a=1}^{\on{dim}\mathfrak{g}} \Delta_{\beta_1,x_av_1}(u_1)\Delta_{\beta_2,x^av_2}(u_2)-\Delta_{x_a\beta_1,v_1}(u_1)\Delta_{x^a\beta_2,v_2}(u_2).
\end{multline}

\begin{defeni}{}\label{defe_univ_class_bethe}
Pick $C \in G$. The (universal) Bethe subalgebra $\ol{{\bf{B}}}(C) \subset \CO(G((z^{-1})))$ is the subalgebra of $\CO(G((z^{-1})))$ generated by the Fourier coefficients of the functions \begin{equation}\label{our_gen_fun_fat_B}
G((z^{-1})) \ni g \mapsto \on{tr}_{V_{\omega_i}} \rho_i(C)\rho_i(g) \in \BC((z^{-1})).
\end{equation}
For $r \in \BZ$, $i \in I$ we denote the coefficient of $z^{-r}$ in (\ref{our_gen_fun_fat_B}) by ${\boldsymbol{\sigma}}_i(C)^{(r)} \in \ol{{\bf{B}}}(C)$. So $\ol{{\bf{B}}}(C)$ is generated by the functions $\{{\boldsymbol{\sigma}}_i(C)^{(r)}$\,\textbar\, $i \in I$, $r \in \BZ$\}. It is convenient to think about (\ref{our_gen_fun_fat_B}) via the generating function
${\boldsymbol{\sigma}}_i(u,C):=\sum_{r \in \BZ} {\boldsymbol{\sigma}}_i(C)^{(r)}u^{-r}.    
$
\end{defeni}

\rem{}
{\em{It is easy to see that the algebra $\ol{{\bf{B}}}(C)$ contains the Fourier coefficients of $g \mapsto \on{tr}_V \rho(C)\rho(g)$ for every finite dimensional representation $\rho\colon G \ra \on{End}(V)$. 
}}
\erem

Let $V$ be a finite dimensional representation of $G$ and pick any basis $e_i \in V$ and the dual basis $e_i^\vee \in V^*$. Using the basis $\{e_i\}_{i=1,\ldots,\on{dim}V}$, we can identify $\on{End}(V)= \on{End}(\BC^n)$.  Let 
\begin{equation*}
{\ol{{\bf{T}}}}^V(u)
\in \on{End}(V) \otimes \CO(G((z^{-1})))[[u,u^{-1}]]=\on{Mat}_{n \times n}\Big(\CO(G((z^{-1})))[[u,u^{-1}]]\Big)
\end{equation*}
be the following element:
\begin{equation*}
{\ol{{\bf{T}}}}^V(u) := 
(\Delta_{e_i^\vee,e_j}(u))_{ij}.
\end{equation*}


We also consider the Casimir element 
$
\Omega:=\sum_{a=1}^{\on{dim}\mathfrak{g}} x_a \otimes x^a \in U(\mathfrak{g})^{\otimes 2}.    
$
The equality (\ref{equal_poiss_in_vectors}) is equivalent to 
\begin{multline}\label{comm_via_T}
\Big\{{\ol{\bf{T}}}^{V_1}_1(u_1),  \ol{{\bf{T}}}_2^{V_2}(u_2)\Big\}=
\\
=\frac{1}{u_2-u_1}\Big((\rho_{V_1} \otimes \rho_{V_2})(\Omega){\bf{\ol{T}}}_1^{V_1}(u_1){\bf{\ol{T}}}_2^{V_2}(u_2) - {\bf{\ol{T}}}_2^{V_2}(u_2){\bf{\ol{T}}}_1^{V_1}(u_1)(\rho_{V_1} \otimes \rho_{V_2})(\Omega)\Big),
\end{multline}
where ${\ol{\bf{T}}}^{V_1}_1(u_1),\, {\ol{\bf{T}}}^{V_2}_2(u_2) \in \on{End}(V_1) \otimes \on{End}(V_2) \otimes \CO(G((z^{-1})))[[u_1,u_1^{-1},u_2,u_2^{-1}]]$ are the images of ${\ol{\bf{T}}}^{V_1}(u_1),\, {\ol{\bf{T}}}^{V_2}(u_2)$ under the embeddings $x \otimes y \mapsto x \otimes 1 \otimes y$, $x \otimes y \mapsto 1 \otimes x \otimes y$, respectively. Indeed, (\ref{equal_poiss_in_vectors}) can be obtained from (\ref{comm_via_T}) by evaluating (\ref{comm_via_T}) on 
$v_1 \in V_1$, $\beta_1 \in V_1^*$, $v_2 \in V_2$, $\beta_2 \in V_2^*$.

\prop{}\label{comm_univ_class_bethe}
The subalgebra $\ol{{\bf{B}}}(C) \subset \CO(G((z^{-1})))$ is Poisson commutative.
\eprop
\prf
We follow \cite[proof of Proposition 4.3]{ir}.
For $i \in I$ set $\ol{{\bf{T}}}^i(u):={\ol{\bf{T}}}^{V_{\omega_i}}(u)$. 
By (\ref{comm_via_T}) we have \begin{equation*}
\Big\{\ol{{\bf{T}}}^i_1(u_1),\ol{{\bf{T}}}^j_2(u_2)\Big\}=\frac{1}{u_2-u_1}\Big((\rho_i \otimes \rho_j)(\Omega){\bf{\ol{T}}}_1^i(u_1){\bf{\ol{T}}}_2^j(u_2) - {\bf{\ol{T}}}_2^j(u_2){\bf{\ol{T}}}_1^i(u_1)(\rho_i \otimes \rho_j)(\Omega) \Big).    
\end{equation*}

We conclude that \begin{multline*}
\Big\{\rho_i(C)_1\ol{{\bf{T}}}^i_1(u_1),\rho_j(C)_2\ol{{\bf{T}}}^j_2(u_2)\Big\}=\\=\frac{\rho_i(C)_1\rho_j(C)_2}{u_2-u_1}\Big((\rho_i \otimes \rho_j)(\Omega){\bf{\ol{T}}}_1^i(u_1){\bf{\ol{T}}}_2^j(u_2) - {\bf{\ol{T}}}_2^j(u_2){\bf{\ol{T}}}_1^i(u_1)(\rho_i \otimes \rho_j)(\Omega) \Big).
\end{multline*}

Note now that $[\rho_i(C) \otimes \rho_j(C),(\rho_i \otimes \rho_j)(\Omega)]=0$. Taking the trace over $V_{\omega_i} \otimes V_{\omega_j}$, using the fact that ${\bf{\ol{T}}}_1^i(u_1)$, ${\bf{\ol{T}}}_2^j(u_2)$ commute (since $\CO(G((z^{-1})))$ is commutative),
we conclude that $\on{tr}_{V_{\omega_i} \otimes V_{\omega_j}}\Big(\Big\{\rho_i(C)_1\ol{{\bf{T}}}^i_1(u_1),\rho_j(C)_2\ol{{\bf{T}}}^j_2(u_2)\Big\}\Big)=0$, thus $\{{{\boldsymbol{\sigma}}}_i(u_1,C),{{\boldsymbol{\sigma}}}_j(u_2,C)\}=0$ as desired.
\epr

Let $G^{\mathrm{reg}} \subset G$ be the subset of regular elements (recall that an element $g \in G$ is called regular if the dimension of the centralizer $Z_G(g)$ is equal to the rank of $\mathfrak{g}$).
\prop{}\label{size_of_univ_class_Bethe}
For $C \in G^{\mathrm{reg}}$ the Fourier coefficients ${\boldsymbol{\sigma}}_i(C)^{(r)}$ of  ${\boldsymbol{\sigma}}_i(u,C)$ are algebraically independent. 
\eprop
\prf
The proof is the same as \cite[proof of Proposition 4.6]{ir}. 
\epr

\subsection{Bethe subalgebras in functions on $\CW_\mu$}{}\label{bethe_in_funct_shift_subsec!}
Let $\mu \in \on{Hom}(\BC^\times,T)=\La$ be an antidominant cocharacter (with respect to a fixed Borel $B \subset G$ containing $T$) and $z^{\mu} \in T((z^{-1}))$ be the corresponding element. Let $G[[z^{-1}]]_1 \subset G[[z^{-1}]]$ be the subgroup consisting of $g(z^{-1})$ such that $g(0)=1 \in G$.
Consider the following (closed) subscheme of $G((z^{-1}))$:
\begin{equation}\label{W_mu_for_antidom}
\CW_\mu:=G[[z^{-1}]]_1 z^\mu  G[[z^{-1}]]_1.    
\end{equation}

\begin{remark}{}\label{W_mu_closed_remark}
{\emph{The scheme $\CW_\mu$ can be defined for arbitrary $\mu \in \La$ as follows: 
\begin{equation}\label{gen_def_W_mu}
\CW_\mu= U[[z^{-1}]]_1 T[[z^{-1}]]_1 z^\mu  U_-[[z^{-1}]]_1,
\end{equation}
where $U$ is the unipotent radical of $B$ and $U_-$ is the unipotent radical of the opposite Borel $B_- \subset G$.
 It is easy to see (see for example~\cite[proof of Theorem A.8]{kpw}) that for antidominant $\mu$ the definition (\ref{gen_def_W_mu}) coincides with (\ref{W_mu_for_antidom}). The fact that $\CW_\mu \subset G((z^{-1}))$ is closed follows from~\cite[Lemma 3.2]{mw}, see also~\cite[Theorem A.8 (a)]{kpw}. To be more precise, it follows from~\cite[Lemma 3.2]{mw} that $\CX_\mu := U((z^{-1}))T[[z^{-1}]]_1 z^\mu U_-((z^{-1})) \subset G((z^{-1}))$ is a closed subfunctor. It remains to note that $\CW_\mu \subset \CX_\mu$ is closed (use that $\CX_\mu \simeq U((z^{-1})) \times T[[z^{-1}]]_1 \times U_-((z^{-1}))$ and $\CW_\mu \simeq U[[z^{-1}]]_1 \times T[[z^{-1}]]_1 \times U_-[[z^{-1}]]_1$).}}
\end{remark}

It follows from \cite[Theorem A.8 (a)]{kpw} that the (closed) embedding $\CW_\mu \subset G((z^{-1}))$ is Poisson (here, it is crucial that $\mu$ is antidominant). 

\rem{}
{\em{Note that we have the decomposition 
\begin{equation*}
G((z^{-1}))=\bigsqcup_{\mu \in \La^-}G[[z^{-1}]] z^\mu G[[z^{-1}]]=\bigsqcup_{\mu \in \La^-}G \cdot \CW_\mu \cdot G,
\end{equation*}
where $\La^- \subset \La$ is the subset of antidominant coweights.}}
\erem

\begin{defeni}{}\label{class_bethe_in_shifted}
Pick $C \in G$ and antidominant $\mu \in \La$. The Bethe subalgebra $\ol{B}_\mu(C) \subset \CO(\CW_\mu)$ is the subalgebra of $\CO(\CW_\mu)$ generated by the Fourier coefficients of the functions \begin{equation}\label{blab!}
\CW_\mu \ni g \mapsto \on{tr}_{V_{\omega_i}} \rho_i(C)\rho_i(g) \in \BC((z^{-1})).
\end{equation}
The coefficient of $z^{-r}$ in (\ref{blab!}) will be denoted by $\sigma_{\mu,i}(C)^{(r)}$. The generating function $\sum_{r \in \BZ}\sigma_{\mu,i}(C)^{(r)} u^{-r}$ is denoted by $\sigma_{\mu,i}(u,C)$. More generally, for every finite dimensional representation $\rho\colon G \ra \on{End}(V)$ we denote by $\sigma_{\mu,V}(C)^{(r)} \in \ol{B}_\mu(C)$ the coefficient of $z^{-r}$ in $g \mapsto \on{tr}_V \rho(C)\rho(g)$.
\end{defeni}

\rem{}
{\em{Note that $\ol{B}_\mu(C)$ is nothing else but the image of $\ol{\bf{B}}(C) \subset \CO(G((z^{-1})))$ under the natural surjection $\CO(G((z^{-1}))) \twoheadrightarrow \CO(\CW_\mu)$ corresponding to the (closed) embedding $\CW_\mu \subset G((z^{-1}))$, $\sigma_{\mu,i}(C)^{(r)}$ is the image of ${\boldsymbol{\sigma}}_i(C)^{(r)}$.}}
\erem

\rem{}
{\em{We follow the notations of \cite{fkprw}. There exists the natural  ``projection'' morphism $\iota_{0,0,\mu}\colon \CW_\mu \twoheadrightarrow \CW_0$ (see \cite[Section 5.9]{fkprw} for details). Morphism $\iota_{0,0,\mu}$ is compatible with Poisson structures. Pullback homomorphism $\iota_{0,0,\mu}^*$ can be quantized to the ``shift'' homomorphism $\iota_{0,0,\mu}\colon Y_0(\mathfrak{g}) \hookrightarrow Y_\mu(\mathfrak{g})$ of (shifted) Yangians of $\mathfrak{g}$, see \cite[Proposition 3.8, Section 5.9, and Theorem 5.15]{fkprw} for details.  Starting from the classical Bethe subalgebra $\ol{B}_0(D) \subset \CO(\CW_0)$, $D \in G$ and taking its pullback $\iota_{0,0,\mu}^*(\ol{B}_0(D))$, we obtain a Poisson commutative subalgebra of $\CO(\CW_\mu)$. It is not true in general that a Bethe subalgebra $\ol{B}_\mu(C) \subset \CO(\CW_\mu)$  is equal to $\iota_{0,0,\mu}^*(\ol{B}_0(D))$ for some $D \in G$ (see Appendix \ref{app_1} for details). On the other hand, one can show that for $C,D \in G^{\mathrm{reg}}$ the Poincar\'e series of $\iota_{0,0,\mu}^*(\ol{B}_0(D))$ are equal to the Poincar\'e series of $\ol{B}_\mu(C)$.}}
\erem

For $\la \in \La$ we set $\la^*:={-}w_0(\la)$, where $w_0 \in W$ is the longest element in the Weyl group of $\mathfrak{g}$. 
\lem{}\label{sigma_zero_neg_degree}
We have $\sigma_{\mu,i}(C)^{(r)}=0$ for $r<\langle \omega_i^*,\mu \rangle$ and  $\sigma_{\mu,i}(C)^{(\langle \omega_i^*,\mu \rangle)}$ is a  positive integer. 
\elem
\prf
This follows from the equality $\CW_\mu=G[[z^{-1}]]_1z^\mu G[[z^{-1}]]_1$.
Indeed, pick $g=g_1 z^\mu g_2 \in G[[z^{-1}]]_1z^\mu G[[z^{-1}]]_1$. 
It is enough to show that if $\la$ is a dominant integral weight and $\rho_\la \colon \mathfrak{g} \ra \on{End}(V_\la)$ is the corresponding irreducible representation, then for every $g \in \CW_\mu$, $\rho_\la(g)$ as a series in $z$ does not have terms of degree greater then $-\langle \la^*,\mu \rangle$ and the term corresponding to $-\langle \la^*,\mu \rangle$ is equal to some positive integer number. Indeed, 
it is clear that $\rho_\la(g_1)$, $\rho_\la(g_2)$ lie in $\on{id}_{V_{\la}}+z^{-1}\on{End}(V_{\la})[[z^{-1}]]$. Recall now that $\mu$ is antidominant so $\rho_{\la}(z^\mu)$ lies in $z^{\langle  w_0(\la),\mu\rangle}\on{End}(V_{\la})[[z^{-1}]]$. It follows that $\rho_\la(g) \in z^{\langle w_0(\la),\mu\rangle}\on{End}(V_\la)[[z^{-1}]]$. It also follows that the coefficient of $z^{\langle w_0(\la),\mu \rangle}$ in $\rho_\la(g)$  is equal to the coefficient of $z^{\langle w_0(\la),\mu \rangle}$ in $\rho_\la(z^{\mu})$.

The coefficient of $z^{\langle w_0(\la),\mu \rangle}$ in $\rho_\la(z^{\mu})$  is the operator $V_{\la} \ra V_{\la}$ that is equal to identity on weight components $V_\la[\nu]$ such that $\langle \nu,\mu \rangle=\langle w_0(\la),\mu \rangle$ and is equal to zero on other components.
The trace of this operator is equal to some positive integer number. 
\epr

Consider the following action $\BC^\times \curvearrowright \CW_{\mu}$: $t \cdot g(z)=g(tz)t^{-\mu}$. We obtain the $\BZ$-grading on $\CO(\CW_\mu)$.

\lem{}\label{degree_delta_our_filtr}
Let $V$ be a finite dimensional representation of $G$ and let $\nu, \nu' \in \on{Hom}(T,\BC^\times)$ be $T$-weights. Pick a vector $v \in V$ of  weight $\nu$ and a covector $\beta \in V^*$ of  weight $\nu'$. Pick $r \in \BZ$ and consider the corresponding function $\Delta_{\beta,v}^{(r)} \in \CO(\CW_\mu)$. Then the degree of this function is equal to $r+ \langle \nu,\mu \rangle$. 
\elem
\prf
We need to compute the action of $t \in \BC^\times$ on the function $\Delta_{\beta,v}^{(r)}$. Recall that for $t \in \BC^\times$, $g \in \CW_\mu$ we have 
\begin{equation*}
t \cdot \Delta_{\beta,v}^{(r)}(g)=\Delta_{\beta,v}^{(r)}(t^{-1} \cdot g)=\Delta_{\beta,v}^{(r)}(g(t^{-1}z)t^{\mu}).
\end{equation*}
So our goal is to compute the  coefficient of $z^{-r}$ in 
\begin{equation*}
\langle \beta, g(t^{-1}z)t^\mu v \rangle = t^{\langle \nu,\mu \rangle}\langle \beta , g(t^{-1}z) v \rangle
\end{equation*}
that is equal to $t^{\langle \nu,\mu \rangle +r}\Delta_{\beta,v}^{(r)}(g)$. So $t \cdot \Delta^{(r)}_{\beta,v}(g)=t^{\langle \nu,\mu \rangle +r}\Delta_{\beta,v}^{(r)}(g)$, hence $\on{deg} \Delta_{\beta,v}^{(r)}(g)=r+\langle \nu,\mu \rangle$.
\epr

Note that the  grading on $\CO(\CW_\mu)$ does not induce the grading on the subalgebra $\ol{B}_\mu(C) \subset \CO(\CW_\mu)$ in general. We always have the induced filtration  $F^\bullet \ol{B}_\mu(C)$ on $\ol{B}_\mu(C)$. 
Note that the graded components of $\CO(\CW_\mu)$ are infinite dimensional in general. As a vector space, the algebra $\ol{B}_\mu(C)$ is spanned by the Fourier coefficients $\sigma_{\mu,V_\la}(C)^{(r)}$ of the functions $g \mapsto \on{tr}_{V_\la}\rho_\la(C)\rho_{\la}(g)$ ($V_\la$ are finite dimensional irreducible representations of $\mathfrak{g}$). The proof of Lemma \ref{sigma_zero_neg_degree} implies that $\sigma_{\mu,V_\la}(C)^{(r)}=0$ for $r<\langle \la^*,\mu \rangle$ and $\sigma_{\mu,V_\la}(C)^{(\langle \la^*,\mu \rangle)} \in \BZ_{>0}$.
It follows from Lemma \ref{degree_delta_our_filtr} that the degree of $\sigma_{\mu,V_\la}(C)^{(r)}$ is at most $r-\langle \la^*,\mu \rangle$. 
In Theorem \ref{size_reg_bethe} below we show that the degree of $\sigma_{\mu,V_\la}(C)^{(r)}$ is precisely $r-\langle \la^*,\mu \rangle$ for $C \in Z_G(\mu)^{\mathrm{reg}}$. As a corollary, we conclude that  $\ol{B}_\mu(C)$ is nonnegatively filtered (with only $\BC$ in zero degree).
So it makes sense to consider Poincar\'e series of $\ol{B}_\mu(C)$ defined as follows:
\begin{equation*}
\on{dim}_q\ol{B}_\mu(C):=\sum_{k \in \BZ} \on{dim}(F^k\ol{B}_\mu(C)/F^{k-1}\ol{B}_\mu(C)) \cdot q^k \in \BZ[[q]].    
\end{equation*}


Using the identification (\ref{gen_def_W_mu}), we have the projection morphism 
\begin{equation*}
\pi\colon \CW_\mu \twoheadrightarrow T[[z^{-1}]]_1z^\mu    
\end{equation*}
that induces the embedding $\pi^*\colon \CO(T[[z^{-1}]]_1z^\mu) \hookrightarrow \CO(\CW_\mu)$.



Set $L:=Z_{G}(\mu)$. Group $L$ is a connected (standard) Levi subgroup of $G$, containing the maximal torus $T$. We denote by $\rho_{L,i}\colon L \ra \on{End}(V_{L,i})$ the irreducible representation of $L$ with the {\emph{lowest weight}} equal to $-\omega_i^*=w_0(\omega_i)$. Note that $V_{L,i}$ can be described as  the irreducible $L$-subrepresentation of $V_{\omega_i}$ generated by the lowest weight vector of $V_{\omega_i}$.


For $C \in L$ and $r>0$ let $\sigma_{L,i}(C)^{(r)}$ be the coefficient of $z^{-r}$ in the function 
\begin{equation*}
L[[z^{-1}]]_1 \ni g \mapsto \on{tr}_{V_{L,i}}\rho_{L,i}(C)\rho_{L,i}(g).
\end{equation*}
We denote by $\ol{B}_L(C) \subset \CO(L[[z^{-1}]]_1)$ the subalgebra of $\CO(L[[z^{-1}]]_1)$ generated by $\sigma_{L,i}(C)^{(r)}$. The algebra $\ol{B}_L(C)$ is nothing else but the standard ``classical'' Bethe subalgebra of $\CO(L[[z^{-1}]]_1)$ (see \cite[Sections 3, 4, 5]{ir}). It follows from 
\cite[Proposition 4.6]{ir} that for $C \in L^{\mathrm{reg}}$ the elements $\{\sigma_{L,i}(C)^{(r)}\,,r>0\}$ are algebraically independent. To be more precise, almost the same proof as the proof of \cite[Proposition 4.6]{ir} works, the only difference is that instead of using \cite[Theorem 3, pg. 119]{st} (see also \cite[Theorem 8.1]{st_art}) we use Proposition \ref{prop_lin_ind_inv} (see Appendix \ref{app}).


\rem{}
{\em{The elements $\sigma_{L,i}(C)^{(r)}$ such that $w_{L,0}w_0(\omega_i)$ is a fundamental weight of $\mathcal{D}(\mathfrak{l}):=[\mathfrak{l},\mathfrak{l}]$ generate the (classical) Bethe subalgebra in $\CO([L,L][[z^{-1}]]_1)$ and the rest of the elements $\sigma_{L,i}(C)^{(r)}$ generate  the (Poisson) center of $\CO(L[[z^{-1}]]_1)$. To see that, one should use the results of the Appendix \ref{app}. We do not provide the details here but refer interested readers to the Appendix \ref{app}.}}
\erem







Consider the closed embedding $L[[z^{-1}]]_1 \subset \CW_\mu$ given by $g \mapsto gz^\mu$.  Note that the restriction of the $\BC^\times$-action on $L[[z^{-1}]]_1$ is just the loop rotation given by $t \cdot g(z)=g(tz)$. 
This embedding induces a surjection $\CO(\CW_\mu) \twoheadrightarrow \CO(L[[z^{-1}]]_1)$ at the level of functions.  Recall that the grading on $\CO(\CW_\mu)$ induces the filtration $F^{\bullet}\ol{B}_\mu(C)$ on $\ol{B}_\mu(C)$. 
Clearly we have natural identifications $\on{gr}\CO(\CW_\mu) \simeq \CO(\CW_\mu)$, $\on{gr}\CO(L[[z^{-1}]]_1) \simeq \CO(L[[z^{-1}]]_1)$. We are now ready to prove Theorem \ref{thm_A}.
\th{}\label{size_reg_bethe}
Assume that $C \in L^{\mathrm{reg}}$. Then the composition 
\begin{equation*}
\on{gr}\ol{B}_\mu(C) \hookrightarrow \CO(\CW_\mu) \twoheadrightarrow \CO(L[[z^{-1}]]_1)
\end{equation*}
induces the identification 
\begin{equation*}
\on{gr}\ol{B}_\mu(C) \iso \ol{B}_L(C).   
\end{equation*}
For $r>\langle \omega_i^*,\mu \rangle$ and $i=1,\ldots,\on{rk}\mathfrak{g}$ the element $\sigma_{\mu,i}(C)^{(r)}$ has degree $r-\langle \omega_i^*,\mu \rangle$ and the identification above  sends $\on{gr}\sigma_{\mu,i}(C)^{(r)}$ to $\sigma_{L,i}(C)^{(r-\langle \omega_i^*,\mu \rangle)}$.   
\eth
\prf
Recall that the functions $\{\sigma_{L,i}(C)^{(r)},\, r>0 \} \subset \CO(L[[z^{-1}]]_1)$ are algebraically independent.
Using Lemma \ref{sigma_zero_neg_degree}, we conclude that we only need to show that 
\begin{equation}\label{gr_restr_sigma}
\on{gr}\Big(\sigma_{\mu,i}(C)^{(r)}|_{L[[z^{-1}]]_1}\Big)=\sigma_{L,i}(C)^{(r-\langle \omega_i^*,\mu \rangle)}.    
\end{equation}

Recall that $V_{L,i}$ is the irreducible representation of $L$ with the {\emph{lowest weight}} equal to $-\omega_i^*=w_0(\omega_i)$. Note that if we decompose $V_{\omega_i}$ as the direct sum of irreducible $L$-modules, then $V_{L,i} \subset V_{\omega_i}$ is precisely the subspace generated by vectors of $V_{\omega_i}$ on which $\mu\colon \BC^\times \ra T$ acts via the multiplication by $t^{-\langle \omega_i^*,\mu\rangle}$, while on the other components it acts via multiplication by some $t^{?}$ with $?<-\langle \omega_i^*,\mu\rangle$ (here we use that $\mu$ is antidominant, cf. proof of Lemma \ref{sigma_zero_neg_degree}). Let $\{v_\nu\}$ be some $T$-weight basis of $V_{\omega_i}$.

For $g \in L[[z^{-1}]]_1$ the value $\sigma_{\mu,i}(C)^{(r)}(gz^\mu)$ is the coefficient of $z^{-r}$ in
\begin{multline*}
\on{tr}_{V_{\omega_i}}\rho_i(C)\rho_i(g)\rho_i(z^\mu)=\sum_{v_\nu}\Delta_{v_\nu^*,v_\nu}(Cgz^{\mu})=
\\
=\sum_{v_\nu: \langle \nu,\mu \rangle=-\langle \omega_i^*,\mu \rangle}\Delta_{v_\nu^*,v_\nu}(Cgz^\mu)+\sum_{v_\nu: \langle \nu,\mu \rangle<-\langle \omega_i^*,\mu \rangle}\Delta_{v_\nu^*,v_\nu}(Cgz^\mu)=\\
=z^{-\langle \omega_i^*,\mu \rangle}\sum_{v_\nu \in V_{L,i}}\Delta_{v_\nu^*,v_\nu}(Cg)+\sum_{v_\nu: \langle \nu,\mu \rangle<-\langle \omega_i^*,\mu \rangle}z^{\langle \nu,\mu \rangle}\Delta_{v_\nu^*,v_\nu}(Cg)
\end{multline*}
that is equal to the coefficient  of  $z^{-r+\langle \omega_i^*,\mu \rangle}$ in 
\begin{equation}\label{main_summand}
\sum_{v_\nu \in V_{L,i}}\Delta_{v_\nu^*,v_\nu}(Cg)=\on{tr}_{V_{L,i}}\rho_i(C)\rho_i(g)
\end{equation}
plus coefficients of $z^{?}$ in some elements $\Delta_{v_\nu^*,v_\nu}(Cg)$ with $?>-r+\langle \omega_i^*,\mu \rangle$. Note now that the degrees (with respect to our filtration) of the terms coming from (\ref{main_summand}) are equal to $r-\langle \omega_i^*,\mu \rangle$ and degrees of the other terms are $<r-\langle \omega_i^*,\mu \rangle$. The equality (\ref{gr_restr_sigma}) follows.
\epr

As a corollary, we conclude that. 

\cor{}\label{cor_size_class_Bethe_shifted}
The elements $\{\sigma_{\mu,i}(C)^{(r)},\, r>\langle \omega_i^*,\mu \rangle\}$ are algebraically independent. The algebra $\ol{B}_\mu(C)$ is a polynomial algebra in the elements $\{\sigma_{\mu,i}(C)^{(r)},\, r>\langle \omega_i^*,\mu \rangle\}$
and the Poincar\'e series of the algebra $\ol{B}_\mu(C)$ coincides with the Poincar\'e series of $\CO(T[[z^{-1}]]_1 z^\mu)$.
\ecor
\prf
Recall that the functions $\{\sigma_{L,i}(C)^{(r)},\, r>0 \} \subset \CO(L[[z^{-1}]]_1)$ are algebraically independent and $\Big(\on{gr}\sigma_{\mu,i}(C)^{(r)}\Big)|_{L[[z^{-1}]]_1}=\sigma_{L,i}(C)^{(r-\langle \omega_i^*,\mu \rangle)}$. Hence, it follows from Theorem \ref{size_reg_bethe} that the elements $\{\sigma_{\mu,i}(C)^{(r)},\, r>\langle \omega_i^*,\mu \rangle\}$ are algebraically independent so the algebra $\ol{B}_\mu(C)$ is a polynomial algebra in them. 
The claim about the Poincar\'e series of $\ol{B}_\mu(C)$ then follows from the definitions.

\epr

\rem{}
{\emph{Note that if $\mu$ is regular then Theorem \ref{size_reg_bethe} tells us that for any $C \in T$ (in particular, for $C=1$) we have $\on{gr} \ol{B}_\mu(C) \simeq \CO(T[[z^{-1}]]_1z^\mu)$ as $\BC^\times$-graded algebras.}}
\erem

\section{Universal Bethe subalgebra in ${\bf{Y}}^{\mathrm{rtt}}(\mathfrak{gl}_n)$}\label{univ_bethe_yang_sect!}
In this section, we assume that $\mathfrak{g}=\mathfrak{gl}_n$. Note that in Section \ref{class-cal_bethe_section} we assumed that $\mathfrak{g}$ is simple, which does not cover the case $\mathfrak{g}=\mathfrak{gl}_n$. So formally the results of Section \ref{class-cal_bethe_section} can not be applied to $\mathfrak{gl}_n$. On the other hand, all of them actually work with only one replacement: we should consider the Poisson algebra $\CO(\on{Mat}_n((z^{-1})))$ instead of $\CO(\on{GL}_n((z^{-1})))$ (see Section \ref{reform_for_gl} below for details). We prefer to deal with $\mathfrak{g}=\mathfrak{gl}_n$ (instead of $\mathfrak{g}=\mathfrak{sl}_n$) since all the objects that appear are simpler in this case, and it is the $\mathfrak{gl}_n$ case not $\mathfrak{sl}_n$ that will be quantized in Sections \ref{univ_bethe_yang_sect!}, \ref{shifted_yangians_section}, \ref{bethe_subalg_in_shifted_section!!!}.
\subsection{Reformulations of the results of Section \ref{class-cal_bethe_section} for $\mathfrak{g}=\mathfrak{gl}_n$}{}\label{reform_for_gl}
For $i,j = 1,\ldots,n$ let $\Delta_{ij}(g)$ be the function on $\mathfrak{gl}_n$ given by $\Delta_{ij}(g)= \langle \epsilon_i^\vee,g\epsilon_j\rangle$. 
For $r \in \BZ$ let $\Delta_{ij}^{(r)} \in \CO(\on{Mat}_n((z^{-1})))$ be the coefficient of $z^{-r}$ in $g(z^{-1}) \mapsto \langle \epsilon_i^\vee, g(z^{-1})\epsilon_j\rangle$. We will denote by the same symbol $\Delta_{ij}^{(r)}$ the restriction of $\Delta_{ij}^{(r)}$ to $\on{GL}_n((z^{-1}))$. We also set $\Delta_{ij}(u):=\sum_{r \in \BZ}\Delta_{ij}^{(r)}u^{-r}$.
The Poisson bracket of $\Delta_{ij}(u)$ on $\CO(\on{GL}_n((z^{-1})))$ is given by the following formula (cf. (\ref{equal_poiss_in_vectors})):
\begin{equation*}
(u_1-u_2)\{\Delta_{ij}(u_1),\Delta_{kl}(u_2)\}=\Delta_{il}(u_1)\Delta_{kj}(u_2)-\Delta_{kj}(u_1)\Delta_{il}(u_2)
\end{equation*}
or equivalently
\begin{equation}\label{Poiss_gl!}
\{\Delta_{ij}^{(p+1)},\Delta_{kl}^{(q)}\}-\{\Delta_{ij}^{(p)},\Delta_{kl}^{(q+1)}\}=\Delta_{il}^{(p)}\Delta_{kj}^{(q)}-\Delta_{kj}^{(p)}\Delta_{il}^{(q)}.
\end{equation}
We see that the bracket above restricts to the Poisson bracket on $\CO(\on{Mat}_n((z^{-1})))$.
We consider $\on{Mat}_n$ as a {\em{monoid}} over $\BC$. We will be interested in finite dimensional representations $V$ of this monoid (i.e. irreducible representations of $\on{GL}_n \subset \on{Mat}_n$ that extend to $\on{Mat}_n$, such representations are called {\em{polynomial}}). Every irreducible representation of $\on{Mat}_n$ can be obtained as a direct summand of the tensor product of representations $\rho_i\colon \on{Mat}_n \ra \on{End}(\La^i\BC^n)$, $i=1,\ldots,n$. We will call $\La^i\BC^n$ {\em{fundamental}} representations of the monoid $\on{Mat}_n$.


The universal Bethe subalgebra of $\CO(\on{Mat}_n((z^{-1})))$ is defined as follows (compare with Definition \ref{defe_univ_class_bethe}).
\begin{defeni}{}
Pick $C \in \on{Mat}_n$. The (universal) Bethe subalgebra $\ol{{\bf{B}}}(C) \subset \CO(\on{Mat}_n((z^{-1})))$ is the subalgebra of $\CO(\on{Mat}_n((z^{-1})))$ generated by the Fourier coefficients of the functions \begin{equation}\label{our_gen_fun_fat_B_gl_case}
\on{Mat}_n((z^{-1})) \ni g \mapsto \on{tr}_{\Lambda^i \BC^n} \rho_i(C)\rho_i(g) \in \BC((z^{-1})).
\end{equation}
For $r \in \BZ$, $i \in \{1,2,\ldots,n\}$ we denote by ${\boldsymbol{\sigma}}_i(C)^{(r)} \in \ol{{\bf{B}}}(C)$ the coefficient of $z^{-r}$ in (\ref{our_gen_fun_fat_B_gl_case}). So $\ol{{\bf{B}}}(C)$ is generated by the functions $\{{\boldsymbol{\sigma}}_i(C)^{(r)}$\,\textbar\, $i \in \{1,2,\ldots,n\}$, $r \in \BZ$\}. It is convenient to think about (\ref{our_gen_fun_fat_B_gl_case}) via the generating function
${\boldsymbol{\sigma}}_i(u,C):=\sum_{r \in \BZ} {\boldsymbol{\sigma}}_i(C)^{(r)}u^{-r}.    
$
\end{defeni}

\rem{}
{\em{It is easy to see that the algebra $\ol{{\bf{B}}}(C)$ contains the Fourier coefficients of $g \mapsto \on{tr}_V \rho(C)\rho(g)$ for every finite dimensional representation $\rho\colon \on{Mat}_n \ra \on{End}(V)$ of the monoid $\on{Mat}_n$.}}
\erem

The following proposition should be compared to Proposition \ref{size_of_univ_class_Bethe}.
\prop{}\label{size_of_univ_class_Bethe_gl}
For $C \in \on{GL}_n^{\mathrm{reg}}$ the Fourier coefficients ${\boldsymbol{\sigma}}_i(C)^{(r)}$ of  ${\boldsymbol{\sigma}}_i(u,C)$ are algebraically independent. 
\eprop
\prf
The proof is almost the same as \cite[proof of Proposition 4.6]{ir}, the only difference is that instead of using \cite[Theorem 3, pg. 119]{st}  we use Proposition \ref{prop_lin_ind_inv}.
\epr

Let $T \subset \on{GL}_n$ be the maximal torus consisting of diagonal matrices. Let $B \subset \on{GL}_n$ be the Borel subgroup consisting of upper triangular matrices. Set $\La:=\on{Hom}(\BC^\times,T)$.
Let $\mu \in \La$ be an antidominant cocharacter (with respect to $B \subset G$) and $z^{\mu} \in T((z^{-1}))$ be the corresponding element. Let $\on{GL}_n[[z^{-1}]]_1 \subset \on{GL}_n[[z^{-1}]]$ be the subgroup consisting of $g(z^{-1})$ such that $g(0)=1 \in \on{GL}_n$.
Consider the following (closed) subscheme of $\on{GL}_n((z^{-1}))$:
\begin{equation}\label{W_mu_for_antidom_gl}
\CW_\mu:=\on{GL}_n[[z^{-1}]]_1 z^\mu  \on{GL}_n[[z^{-1}]]_1.   
\end{equation}

\rem{}
{\em{Note that $\CW_\mu$ is closed in both $\on{GL}_n((z^{-1}))$ and $\on{Mat}_n((z^{-1}))$, the argument is similar to the one in Remark \ref{W_mu_closed_remark}. 
}}
\erem

\begin{defeni}{}\label{class_bethe_in_shifted_gl}
Pick $C \in \on{GL}_n$ and antidominant $\mu \in \La$. The classical Bethe subalgebra $\ol{B}_\mu(C) \subset \CO(\CW_\mu)$ is the subalgebra of $\CO(\CW_\mu)$ generated by the Fourier coefficients of the functions \begin{equation}\label{blump!}
\CW_\mu \ni g \mapsto \on{tr}_{\Lambda^i\BC^n} \rho_i(C)\rho_i(g) \in \BC((z^{-1})).
\end{equation}
The coefficient of $z^{-r}$ in (\ref{blump!}) will be denoted by $\sigma_{\mu,i}(C)^{(r)}$. The generating function $\sum_{r \in \BZ}\sigma_{\mu,i}(C)^{(r)} u^{-r}$ is denoted by $\sigma_{\mu,i}(u,C)$. More generally, for every finite dimensional representation $\rho\colon \on{Mat}_n \ra \on{End}(V)$ we denote by $\sigma_{\mu,V}(C)^{(r)} \in \ol{B}_\mu(C)$ the coefficient of $z^{-r}$ in $g \mapsto \on{tr}_V \rho(C)\rho(g)$.
\end{defeni}

We identify naturally $\on{Lie}T \simeq \BC^n$ and denote by $\epsilon_1,\ldots,\epsilon_n$ the standard basis of $\BC^n$. Recall that $\epsilon_1^{\vee},\ldots,\epsilon_n^{\vee} \in (\BC^n)^*$ is the dual basis. For $k=1,\ldots,n$ set $\omega_k^*:=-\epsilon_n^{\vee}-\ldots-\epsilon_{n-k+1}^{\vee}$. Set $L:=Z_{\on{GL}_n}(\mu)$. The following theorem holds.
\th{}\label{size_reg_bethe_gl}
Assume that $C \in L^{\mathrm{reg}}$. Then the composition 
\begin{equation*}
\on{gr}\ol{B}_\mu(C) \hookrightarrow \CO(\CW_\mu) \twoheadrightarrow \CO(L[[z^{-1}]]_1)
\end{equation*}
induces the identification 
\begin{equation*}
\on{gr}\ol{B}_\mu(C) \iso \ol{B}_L(C).
\end{equation*}
For $r > \langle \omega_i^*,\mu \rangle$ and $i=1,\ldots,n$ the element $\sigma_{\mu,i}(C)^{(r)}$ has degree $r-\langle \omega_i^*,\mu \rangle$ and the identification above  sends $\on{gr}\sigma_{\mu,i}(C)^{(r)}$ to $\sigma_{L,i}(C)^{(r-\langle \omega_i^*,\mu \rangle)}$.   
\eth
\prf
The proof is the same as the one of Theorem \ref{size_reg_bethe}.
\epr

\cor{}\label{cor_size_class_Bethe_shifted_gl}
The elements $\{\sigma_{\mu,i}(C)^{(r)},\, r>\langle \omega_i^*,\mu \rangle\}$ are algebraically independent. The algebra $\ol{B}_\mu(C)$ is a polynomial algebra in the elements $\{\sigma_{\mu,i}(C)^{(r)},\, r>\langle \omega_i^*,\mu \rangle\}$
and the Poincar\'e series of the algebra $\ol{B}_\mu(C)$ coincides with the Poincar\'e series of $\CO(T[[z^{-1}]]_1 z^\mu)$.
\ecor
\prf
The proof is the same as the one of Corollary \ref{cor_size_class_Bethe_shifted}.
\epr

\subsection{The ``Yangian'' quantization of $\CO(\on{Mat}_n((z^{-1})))$}{}
Consider the rational $R$-matrix $R(u)$ for $\mathfrak{g}=\mathfrak{gl}_n$:
\begin{equation*}
R(u)=-u - P,\, P=\sum_{i,j}E_{ij} \otimes E_{ji},   
\end{equation*}
where $E_{ij} \in \mathfrak{gl}_n$ are the elementary matrices.

\rem{}\label{P_vs_Casimir}
{\em{Note that $P \in U(\mathfrak{gl}_n)^{\otimes 2}$ is nothing else but the Casimir element $\Omega$ corresponding to $\mathfrak{gl}_n$.}}
\erem

\begin{defeni}{}\label{defe_RTT_Yang_univ!}
Let ${\bf{Y}}^{\mathrm{rtt}}(\mathfrak{gl}_n)^{\mathrm{pol}}$ be the associative $\BC$-algebra generated by $\{t_{ij}^{(r)}\}^{r \in \BZ}_{1 \leqslant i,j \leqslant n}$ subject to the following  family of relations:
\begin{equation}\label{RTT_univ!}
R(u_1-u_2){\bf{T}}_1(u_1){\bf{T}}_2(u_2)={\bf{T}}_2(u_2){\bf{T}}_1(u_1)R(u_1-u_2),    
\end{equation}
where ${\bf{T}}(u) \in {\bf{Y}}^{\mathrm{rtt}}(\mathfrak{gl}_n)^{\mathrm{pol}}[[u,u^{-1}]] \otimes \on{End}(\BC^n)$ is defined via 
\begin{equation*}
{\bf{T}}(u)=(t_{ij}(u))_{ij}~\text{with}~t_{ij}(u)=\sum_{r \in \BZ} t_{ij}^{(r)}u^{-r}.    
\end{equation*}

\rem{}\label{rem_rtt_via_poiss}
{\emph{Note that  (\ref{RTT_univ!}) is equivalent to the equality \begin{equation}\label{comm_poiss_rtt}
(u_2-u_1)[{\bf{T}}_1(u_1),{\bf{T}}_2(u_2)]=P{\bf{T}}_1(u_1){\bf{T}}_2(u_2)-{\bf{T}}_2(u_2){\bf{T}}_1(u_1)P.   
\end{equation}
This equation is a quantization of (\ref{comm_via_T}). Explicitly, (\ref{comm_poiss_rtt}) is equivalent to the  following relations:
\begin{equation}\label{def_rel_thick}
[t_{ij}^{(p+1)},t_{kl}^{(q)}]-[t_{ij}^{(p)},t_{kl}^{(q+1)}]=t_{kj}^{(q)}t_{il}^{(p)}-t_{kj}^{(p)}t_{il}^{(q)},\, p,q \in \BZ.    
\end{equation}
}}
\erem

For every $N \in \BZ_{\geqslant 0}$ consider the two-sided ideal $I_N \subset {\bf{Y}}^{\mathrm{rtt}}(\mathfrak{gl}_n)^{\mathrm{pol}}$ generated by the elements 
$
\{t_{ij}^{(-r)}\,|\, r > N\}.
$
We set 
\begin{equation*}
{\bf{Y}}^{\mathrm{rtt}}(\mathfrak{gl}_n)_{N}:={\bf{Y}}^{\mathrm{rtt}}(\mathfrak{gl}_n)^{\mathrm{pol}}/I_N.    
\end{equation*}
Let ${\bf{Y}}^{\mathrm{rtt}}(\mathfrak{gl}_n)$ be the corresponding completion of ${\bf{Y}}^{\mathrm{rtt}}_\mu(\mathfrak{gl}_n)^{\mathrm{pol}}$:
\begin{equation*}
{\bf{Y}}^{\mathrm{rtt}}(\mathfrak{gl}_n):=\underset{\longleftarrow}{\on{lim}}\,{\bf{Y}}^{\mathrm{rtt}}(\mathfrak{gl}_n)_N.
\end{equation*}
\end{defeni}


\lem{}\label{shift_identif}
We have an isomorphism of algebras $Y^{\mathrm{rtt}}(\mathfrak{gl}_n)_N \iso Y^{\mathrm{rtt}}(\mathfrak{gl}_n)_0$ given by $t_{ij}^{(r)} \mapsto t_{ij}^{(r+N)}$.
\elem
\prf
Follows from the definitions.
\epr

\lem{}\label{coprod_Y_0}
The assignment $\Delta({\bf{T}}(u))={\bf{T}}(u) \otimes {\bf{T}}(u)$ or equivalently 
\begin{equation*}
\Delta(t_{ij}^{(r)})=\sum_{k=1}^n \sum_{p+q=r} t_{ik}^{(p)} \otimes t_{kj}^{(q)}    
\end{equation*}
extends to the homomorphism of algebras $\Delta\colon {\bf{Y}}^{\mathrm{rtt}}(\mathfrak{gl}_n)_0 \ra {\bf{Y}}^{\mathrm{rtt}}(\mathfrak{gl}_n)_0 \otimes {\bf{Y}}^{\mathrm{rtt}}(\mathfrak{gl}_n)_0$.
\elem
\prf
Follows from the definitions.
\epr

As before, let  $\epsilon_1,\ldots,\epsilon_n$ be the standard basis of $\BC^n$,  let $\epsilon_1^{\vee},\ldots,\epsilon_n^{\vee} \in (\BC^n)^*$ be the dual basis, and recall $\Delta_{ij}=\Delta_{\epsilon_i^\vee,\epsilon_j} \in \CO(\on{Mat}_n)$.

\lem{}\label{coev_gl_n}
We have a (surjective) homomorphism of algebras ${\bf{Y}}^{\mathrm{rtt}}(\mathfrak{gl}_n)_0 \twoheadrightarrow \CO(\on{Mat}_n)$ given by $t_{ij}^{(r)} \mapsto \delta_{0,r}\Delta_{ij}$.
\elem
\prf
Recall that the algebra ${\bf{Y}}^{\mathrm{rtt}}(\mathfrak{gl}_n)_0$ is generated by $t_{ij}^{(r)}$, $r \in \BZ$ subject to relations 
\begin{equation}\label{eq_def_1}
[t_{ij}^{(p+1)},t_{kl}^{(q)}]-[t_{ij}^{(p)},t_{kl}^{(q+1)}]=t_{kj}^{(q)}t_{il}^{(p)}-t_{kj}^{(p)}t_{il}^{(q)},\, p,q \in \BZ,    
\end{equation}
\begin{equation*}
t_{ij}^{(r)}=0,~\text{for}~r<0.    
\end{equation*}
Note now that the image in $\CO(\on{Mat}_n)$ of the LHS of the equality (\ref{eq_def_1}) is equal to zero since the algebra $\CO(\on{Mat}_n)$ is commutative. It remains to note that the image of the RHS is 
\begin{equation*}
\delta_{0,q}\Delta_{kj}\delta_{0,p}\Delta_{il}-\delta_{0,p}\Delta_{kj}\delta_{0,q}\Delta_{il}=0.
\end{equation*}
\epr

\lem{}\label{lem_hom_fat_to_thick_yangian}
We have a (surjective) homomorphism ${\bf{Y}}^{\mathrm{rtt}}(\mathfrak{gl}_n)_{0} \twoheadrightarrow Y(\mathfrak{gl}_n)$ given by 
\begin{equation*}
t_{ij}^{(0)} \mapsto \delta_{ij},\, t_{ij}^{(r)} \mapsto t_{ij}^{(r)}~\text{for}~r>0.    
\end{equation*}
\elem
\prf
Follows from the definitions.
\epr

Composing the homomorphisms defined in Lemmas \ref{shift_identif}, \ref{coprod_Y_0}, \ref{coev_gl_n}, \ref{lem_hom_fat_to_thick_yangian} we obtain the homomorphism 
\begin{equation*}
{\bf{Y}}^{\mathrm{rtt}}(\mathfrak{gl}_n)_{N} \iso  {\bf{Y}}^{\mathrm{rtt}}(\mathfrak{gl}_n)_{0} \xrightarrow{\Delta}   {\bf{Y}}^{\mathrm{rtt}}(\mathfrak{gl}_n)_{0} \otimes {\bf{Y}}^{\mathrm{rtt}}(\mathfrak{gl}_n)_{0} \ra \CO(\on{Mat}_n) \otimes Y(\mathfrak{gl}_n)
\end{equation*}
given by: 
\begin{equation*}
t_{ij}^{(r)} \mapsto t_{ij}^{(r+N)} \mapsto \sum_{k=1}^{n}\sum_{p+q=r+N,\, p,q \geqslant 0}t_{ik}^{(p)} \otimes t_{kj}^{(q)}\mapsto \sum_{k=1}^{n}\Delta_{ik} \otimes t^{(r+N)}_{kj}~\text{for}~r>-N, 
\end{equation*}
\begin{equation*}
t_{ij}^{(-N)} \mapsto \sum_{k=1}^n \Delta_{ik} \otimes \delta_{kj}=\Delta_{ij} \otimes 1.
\end{equation*}
We denote this homomorphism by $\Psi_N$.


The algebras ${\bf{Y}}^{\mathrm{rtt}}(\mathfrak{gl}_n)_N$ are naturally filtered by placing $t_{ij}^{(r)}$ in degree $r$. 

\prop{}(PBW theorem for ${\bf{Y}}^{\mathrm{rtt}}(\mathfrak{gl}_n)_{N}$)\label{gen_quot_yang_univ}
 The ordered products of the elements of the set $\{t_{ij}^{(r)},\, r \geqslant -N,\, 1 \leqslant i,j \leqslant n\}$ form a $\mathbb{C}$-basis of  ${\bf{Y}}^{\mathrm{rtt}}(\mathfrak{gl}_n)_{N}$. The algebra $\on{gr}{\bf{Y}}^{\mathrm{rtt}}(\mathfrak{gl}_n)_{N}$ is a polynomial algebra generated by the elements $\{\on{gr}t_{ij}^{(r)},\, r \geqslant -N\}$.
\eprop
\prf
Let us first of all recall that by Lemma \ref{shift_identif} ${\bf{Y}}^{\mathrm{rtt}}(\mathfrak{gl}_n)_{N} \simeq {\bf{Y}}^{\mathrm{rtt}}(\mathfrak{gl}_n)_{0}$, hence, we can assume that $N=0$. It follows from the definitions that $\on{gr}{\bf{Y}}^{\mathrm{rtt}}(\mathfrak{gl}_n)_{0}$ is commutative and the polynomials in $\{t_{ij}^{(r)},\, r \geqslant 0,\, 1 \leqslant i,j \leqslant n\}$ indeed span the algebra ${\bf{Y}}^{\mathrm{rtt}}(\mathfrak{gl}_n)_{0}$. 
 It remains to show that the elements $\{\on{gr}t_{ij}^{(r)},\, r \geqslant 0,\, 1 \leqslant i,j \leqslant n\}$ are algebraically independent.
The algebra $Y(\mathfrak{gl}_n)$ is filtered via $\on{deg}t_{ij}^{(r)}=r$ and this filtration induces the filtration on $\CO(\on{Mat}_n) \otimes Y(\mathfrak{gl}_n)$  (placing $\CO(\on{Mat}_n)$ in degree zero). 
The homomorphism 
\begin{equation*}
\Psi_0\colon {\bf{Y}}^{\mathrm{rtt}}(\mathfrak{gl}_n)_0 \ra \CO(\on{Mat}_n) \otimes Y(\mathfrak{gl}_n) \end{equation*}
 is filtered.  We have
\begin{equation*}
\on{gr}\Big( \CO(\on{Mat}_n) \otimes Y(\mathfrak{gl}_n)\Big)=\CO(\on{Mat}_n) \otimes \CO(1+z^{-1}\on{Mat}_n[[z^{-1}]])
\end{equation*}
and from the PBW theorems for $\CO(\on{Mat}_n)$, $Y(\mathfrak{gl}_n)$ we conclude that $\CO(\on{Mat}_n) \otimes Y(\mathfrak{gl}_n)$ has the following basis:
\begin{equation*}
\{(\Delta_{k_1l_1} \cdot \ldots \cdot \Delta_{k_pl_p}) \otimes (t_{i_1j_1}^{(r_1)} \cdot \ldots \cdot t_{i_sj_s}^{(r_s)})\,|\, \Delta_{k_1l_1}< \ldots < \Delta_{k_pl_p},\, t_{i_1j_1}^{(r_1)} < \ldots < t_{i_sj_s}^{(r_s)},\, p,s \in \BZ_{\geqslant 0}\},
\end{equation*}
where $<$ are any orders on the sets $\{\Delta_{kl}\}_{1 \leqslant k,l \leqslant n}$, $\{t_{ij}^{(r)}\}_{1 \leqslant i,j \leqslant n,\, 1 \leqslant r}$.

It remains to note that the composition (where the first morphism is the natural embedding and $m$ is the multiplication morphism)
\begin{equation*}
\on{Mat}_n \times (1+z^{-1}\on{Mat}_n[[z^{-1}]]) \ra \on{Mat}_n[[z^{-1}]] \times \on{Mat}_n[[z^{-1}]] \xrightarrow{m}  \on{Mat}_n[[z^{-1}]]   
\end{equation*}
contains $\on{GL}_n[[z^{-1}]]$ in the image so is dominant, hence,  induces the embedding at the level of functions. Since the functions $\{\Delta_{ij}^{(r)},\,r \geqslant 0,\, 1 \leqslant i,j \leqslant n\}$ are algebraically independent functions on $\on{Mat}_n[[z^{-1}]]$ with pullbacks exactly  $\on{gr}\Big(\sum_{k=1}^{n}\Delta_{ik} \otimes t^{(r)}_{kj}\Big)$, this implies that the elements $\on{gr}\Big(\sum_{k=1}^{n}\Delta_{ik} \otimes t^{(r)}_{kj}\Big)$ are algebraically independent. It follows that the elements $\{\on{gr}t_{ij}^{(r)} ,\, r \geqslant 0,\, 1 \leqslant i,j \leqslant n\}$ are algebraically independent (use that $(\on{gr}\Psi_0)(\on{gr}t_{ij}^{(r)})=\on{gr}\Big(\sum_{k=1}^{n}\Delta_{ik} \otimes t^{(r)}_{kj}\Big)$).
\epr


\prop{}\label{iso_poiss_gr_yang_gl_univ}
We have an isomorphism of Poisson algebras 
\begin{equation*}
\on{gr}{\bf{Y}}^{\mathrm{rtt}}(\mathfrak{gl}_n) \iso \CO(\on{Mat}_n((z^{-1})))    
\end{equation*}
given by $\on{gr}t_{ij}^{(r)} \mapsto \Delta_{ij}^{(r)}$, where $\on{gr}{\bf{Y}}^{\mathrm{rtt}}(\mathfrak{gl}_n):=\underset{\longleftarrow}{\on{lim}}\,\on{gr}{\bf{Y}}^{\mathrm{rtt}}_\mu(\mathfrak{gl}_n)_N$.
\eprop
\prf
Let us first of all show that the map $\on{gr}t_{ij}^{(r)} \mapsto \Delta_{ij}^{(r)}$ defines an isomorphism of algebras $\on{gr}{\bf{Y}}^{\mathrm{rtt}}_\mu(\mathfrak{gl}_n) \iso \CO(\on{Mat}_n((z^{-1})))$. Recall that 
\begin{equation*}
{\bf{Y}}^{\mathrm{rtt}}(\mathfrak{gl}_n)=\underset{\longleftarrow}{\on{lim}}\,{\bf{Y}}^{\mathrm{rtt}}(\mathfrak{gl}_n)_N,\, \CO(\on{Mat}_n((z^{-1}))) = \underset{\longleftarrow}{\on{lim}}\,\CO(z^{N}\on{Mat}_n[[z^{-1}]])
\end{equation*}
so it remains to show that the map $\on{gr}t_{ij}^{(r)} \mapsto \Delta_{ij}^{(r)}$ defines an isomorphism of algebras $\on{gr}{\bf{Y}}^{\mathrm{rtt}}(\mathfrak{gl}_n)_N \iso \CO(z^{N}\on{Mat}_n[[z^{-1}]])$. This is clear since both of these algebras are polynomial algebras in the corresponding generators (here we use Proposition \ref{gen_quot_yang_univ}). The fact that the isomorphism is Poisson follows from the definitions (see Equations (\ref{comm_via_T}), (\ref{RTT_univ!}) and also Remarks \ref{P_vs_Casimir}, \ref{rem_rtt_via_poiss}). 
\epr


It follows from Proposition \ref{iso_poiss_gr_yang_gl_univ} that the algebra ${\bf{Y}}^{\mathrm{rtt}}(\mathfrak{gl}_n)$ is a (filtered) quantization of the Poisson algebra $\CO(\on{Mat}_n((z^{-1})))$. Recall now that in Section \ref{reform_for_gl} (see also Section \ref{class-cal_bethe_section}) we considered closed Poisson subschemes $\CW_\mu \subset \on{Mat}_n((z^{-1}))$. We also considered the action of $\BC^\times$ on $\CW_\mu$ given by $g(z) \mapsto g(tz)t^{-\mu}$ for $g(z) \in \CW_\mu$. Note that more generally for every $\mu_1, \mu_2 \in \La$ such that $\mu_1+\mu_2=\mu$ we have a $\BC^\times$-action on $\CW_\mu$ given by $g(z) \mapsto t^{-\mu_1}g(tz)t^{-\mu_2}$, cf. \cite[Section 5.9]{fkprw}.  This action can be lifted to the action on the whole $\on{Mat}_n((z^{-1}))$. Moreover, note that we have an action of $\BC^\times \times T \times T$ on $\on{Mat}_n((z^{-1}))$ given by 
\begin{equation*}
(t,s_1,s_2) \cdot g(z) = s_1^{-1}g(tz)s_2^{-1},\, (t,s_1,s_2) \in \BC^\times \times T \times T.    
\end{equation*}
Then the former $\BC^\times$-action above just corresponds to the cocharacter $\BC^\times \ra \BC^\times \times T \times T$ given by $t \mapsto (t,t^{\mu_1},t^{\mu_2})$. Note also that by Proposition \ref{iso_poiss_gr_yang_gl_univ}  the action of $\BC^\times \times \{1\} \times \{1\} \subset \BC^\times \times T \times T$ induces the grading on $\CO(\on{Mat}_n((z^{-1})))$ that lifts to the {\em{filtration}} on ${\bf{Y}}^{\mathrm{rtt}}(\mathfrak{gl}_n)$ given by $\on{deg}t_{ij}^{(r)}=r$. We now claim that:

\lem{}\label{La_La_grad_univ} The action of $T \times T$ on $\on{Mat}_n((z^{-1}))$ induces the $\La \times \La$-grading that lifts to the {\emph{$\La \times \La$-grading}} on ${\bf{Y}}^{\mathrm{rtt}}(\mathfrak{gl}_n)$ given by the following formula:
\begin{equation*}
\on{deg}t_{ij}^{(r)}=((0,\ldots,0,\underset{i}{1},0,\ldots,0),(0,\ldots,0,\underset{j}{1},0,\ldots,0)) \in \La \times \La,     
\end{equation*}
where we naturally identify $\La \simeq \BZ^n$.
\elem
\prf
Recall that ${\bf{Y}}^{\mathrm{rtt}}(\mathfrak{gl}_n)$ is the inverse limit of the algebras ${\bf{Y}}^{\mathrm{rtt}}(\mathfrak{gl}_n)_N$. The algebra ${\bf{Y}}^{\mathrm{rtt}}(\mathfrak{gl}_n)_N$ is generated by $t_{ij}^{(r)}$, $r \in \BZ$, subject to relations:
\begin{equation*}
[t_{ij}^{(p+1)},t_{kl}^{(q)}]-[t_{ij}^{(p)},t_{kl}^{(q+1)}]=t_{kj}^{(q)}t_{il}^{(p)}-t_{kj}^{(p)}t_{il}^{(q)},\, p,q \in \BZ,    
\end{equation*}
\begin{equation*}
t_{ij}^{(r)}=0,~\text{for}~r<-N.    
\end{equation*}
All of these relations are homogeneous (with respect to the above $\La \times \La$-grading). It is also clear that the $\La \times \La$-grading above is compatible with the filtration on ${\bf{Y}}^{\mathrm{rtt}}(\mathfrak{gl}_n)_N$  given by $\on{deg}t_{ij}^{(r)}=r$. It remains to note that the degree of $\Delta_{ij}^{(r)} \in \CO(\on{Mat}_n((z^{-1})))$ with respect to the $\La \times \La$-grading induced by the $T \times T$-action on $\on{Mat}_n((z^{-1}))$ is equal to $((0,\ldots,0,\underset{i}{1},0,\ldots,0),(0,\ldots,0,\underset{j}{1},0,\ldots,0))$. 
 \epr

Using Lemma \ref{La_La_grad_univ} we can now define the ``$\mu$-twisted'' filtration on ${\bf{Y}}^{\mathrm{rtt}}(\mathfrak{gl}_n)$ as follows: for $\mu=(d_1,\ldots,d_n) \in \BZ^n=\La$ we set
\begin{equation}\label{filtr_twist_univ}
\on{deg}_\mu(t_{ij}^{(r)}):=r+d_j.   
\end{equation}

\rem{}
{\em{Note that more generally for every $\mu_1,\mu_2 \in \La$, we can define the filtration $\on{deg}_{\mu_1,\mu_2}(t_{ij}^{(r)}):=r+a_i+b_j$, where $\mu_1=(a_1,\ldots,a_n) \in \BZ^n$, $\mu_2=(b_1,\ldots,b_n) \in \BZ^n$. The filtration $\on{deg}_\mu$ that we consider corresponds to taking $\mu_1=0$, $\mu_2=\mu$ and at the classical level  to the $\BC^\times$-action on $\on{Mat}_n((z^{-1}))$ via $g(z) \mapsto g(tz)t^{-\mu}$ (recall that the $\BC^\times$-action on $\on{Mat}_n((z^{-1}))$ sends a function $f(g)$ to the function $g \mapsto f(t^{-1}g)$).}}
\erem

In Section \ref{shifted_yangians_section}, we will discuss certain quotients of ${\bf{Y}}^{\mathrm{rtt}}(\mathfrak{gl}_n)$ called shifted Yangians. These algebras quantize $\CO(\CW_\mu)$ and are equipped with the natural filtration quantizing the $\BC^\times$-action on $\CW_\mu$ given by $g(z) \mapsto g(tz)t^{-\mu}$. This filtration is compatible with the ``$\mu$-twisted'' filtration (\ref{filtr_twist_univ}) on ${\bf{Y}}^{\mathrm{rtt}}(\mathfrak{gl}_n)$ (see Lemma \ref{gr_t_i_j!} below).

Let us finish this section with the following proposition that is just the generalization of Proposition \ref{iso_poiss_gr_yang_gl_univ} to the case of ``$\mu$-twisted'' filtration (\ref{filtr_twist_univ}).

\prop{}\label{iso_poiss_gr_yang_gl_univ_twisted}
For every $\mu \in \La$ we have an isomorphism of graded Poisson algebras 
\begin{equation*}
\on{gr}_\mu{\bf{Y}}^{\mathrm{rtt}}(\mathfrak{gl}_n) \iso \CO(\on{Mat}_n((z^{-1})))  
\end{equation*}
given by $\on{gr}_\mu t_{ij}^{(r)} \mapsto \Delta_{ij}^{(r)}$, where $\on{gr}_\mu{\bf{Y}}^{\mathrm{rtt}}(\mathfrak{gl}_n):=\underset{\longleftarrow}{\on{lim}}\,\on{gr}_\mu{\bf{Y}}^{\mathrm{rtt}}_\mu(\mathfrak{gl}_n)_N$, and the grading on $\CO(\on{Mat}_n((z^{-1})))$ corresponds to the $\BC^\times$-action on $\on{Mat}_n((z^{-1}))$ given by $g(z) \mapsto g(tz)t^{-\mu}$.
\eprop
\prf
Follows from Proposition \ref{iso_poiss_gr_yang_gl_univ} and Lemma \ref{La_La_grad_univ} (recall that $\on{gr}_{\mu} t_{il}^{(p)} \on{gr}_\mu t_{kj}^{(q)}=\on{gr}_\mu t_{kj}^{(q)}\on{gr}_{\mu} t_{il}^{(p)}$). Let us explain why the isomorphism is Poisson. Note that the term $[t_{ij}^{(p+1)},t_{kl}^{(q)}]-[t_{ij}^{(p)},t_{kl}^{(q+1)}]$ in (\ref{def_rel_thick}) has degree $p+q+d_j+d_l+1$, and the degree of $t_{kj}^{(q)}t_{il}^{(p)}-t_{kj}^{(p)}t_{il}^{(q)}$ is equal to $p+q+d_j+d_l$.
Passing to the associated graded it then follows from (\ref{def_rel_thick}) that the Poisson structure on $\on{gr}_\mu {\bf{Y}}^{\mathrm{rtt}}(\mathfrak{gl}_n)$ is given by 
\begin{equation}\label{rel_gr_Poiss}
\{\on{gr}_\mu t_{ij}^{(p+1)},\on{gr}_\mu t_{kl}^{(q)}\}-\{\on{gr}_\mu t_{ij}^{(p)},\on{gr}_\mu t_{kl}^{(q+1)}\}=\on{gr}_\mu t_{il}^{(p)}\on{gr}_\mu t_{kj}^{(q)}-\on{gr}_\mu t_{kj}^{(p)}\on{gr}_\mu t_{il}^{(q)}.
\end{equation}
It remains to note that (\ref{rel_gr_Poiss}) clearly coincides with (\ref{Poiss_gl!}) after the identification $\on{gr}_\mu t_{ij}^{(r)} \mapsto \Delta_{ij}^{(r)}$ thus implying the compatibility of Poisson structures.   
\epr



\subsection{Universal Bethe subalgebras ${\bf{B}}(C) \subset {\bf{Y}}^{\mathrm{rtt}}(\mathfrak{gl}_n)$}{}\label{univ_bethe_in_yang} For every $k=1,\ldots,n$ we denote by $A_k \in \on{End}(\BC^n)^{\otimes k}$ the antisymmetrization  map normalized so that $A_k^2=A_k$. For $C \in \on{GL}_n$ we set
\begin{equation*}
\boldsymbol{\tau}_{k}(u,C)_N:=\on{tr}_{{(\BC^n)^{\otimes k}}} A_k C_1 \ldots C_k {\bf{T}}_1(u) \ldots {\bf{T}}_k(u-k+1) \in {\bf{Y}}^{\mathrm{rtt}}(\mathfrak{gl}_n)_N((u^{-1})),  
\end{equation*}
where $C_i \in \on{End}(\BC^n)^{\otimes k}$ is the matrix $1 \otimes \ldots \otimes 1 \otimes \underset{i}{C} \otimes 1 \otimes \ldots \otimes 1$.

\begin{defeni}{}\label{defe_bethe_univ}
We denote by ${\bf{B}}(C)_N \subset {\bf{Y}}^{\mathrm{rtt}}(\mathfrak{gl}_n)_N$ the subalgebra generated by the coefficients of $\boldsymbol{\tau}_{k}(u,C)_N$, $k=1,\ldots,n$. Let ${\bf{B}}(C)$ be the subalgebra of  $\underset{\longleftarrow}{\on{lim}}\,{\bf{B}}(C)_N \subset {\bf{Y}}^{\mathrm{rtt}}(\mathfrak{gl}_n)$ generated by the coefficients of the elements $\boldsymbol{\tau}_{k}(u,C)=\underset{\longleftarrow}{\on{lim}}\,\boldsymbol{\tau}_{k}(u,C)_N$. We call ${\bf{B}}(C)$ the {\emph{universal Bethe subalgebra}} of ${\bf{Y}}^{\mathrm{rtt}}(\mathfrak{gl}_n)$.
\end{defeni}

Let us show that the algebra ${\bf{B}}(C)$ is commutative. It is enough to show that the algebras ${\bf{B}}(C)_N$ are commutative. This is a direct corollary of the following proposition.
\prop{}\label{univ_bethe_are_comm}For any $N$ the coefficients of the series $\{\boldsymbol{\tau}_{k}(u,C)_N\}_{k=1}^{n}$ pairwise commute.
\eprop
\prf
The proof is precisely the same as the one of the fact that $[\hat{B}_k(u),\hat{B}_l(v)]=0$ in \cite[Proposition 1.2]{no}, as it utilized only the RTT relation (\ref{RTT_univ!}) and the Yang-Baxter relation $R_{12}(u)R_{13}(u+v)R_{23}(v)=R_{23}(v)R_{13}(u+v)R_{12}(u)$ satisfied by $R(u)$.
\epr

 
\begin{Ex}
For $\mathfrak{g}=\mathfrak{gl}_2$ and $C=\on{diag}(c_1,c_2)$
we get:
\begin{equation*}
\boldsymbol{\tau}_{1}(u,C)_N=\on{tr}C{\bf{T}}(u)=c_1t_{11}(u)+c_2t_{22}(u),   
\end{equation*}
\begin{multline*}
\boldsymbol{\tau}_{2}(u,C)_N=\on{tr} A_2 C_1C_2 {\bf{T}}_1(u){\bf{T}}_2(u-1)=\\
=\frac{1}{2}c_1c_2(t_{11}(u)t_{22}(u-1)+t_{22}(u)t_{11}(u-1)-t_{12}(u)t_{21}(u-1)-t_{21}(u)t_{12}(u-1)).    
\end{multline*}

\end{Ex}

Recall that the algebra ${\bf{Y}}^{\mathrm{rtt}}(\mathfrak{gl}_n)$ is filtered via $\on{deg}t_{ij}^{(r)}=r$. We obtain a filtration on ${\bf{B}}(C)$.
\prop{}\label{gen_of_fat_B_prop}
Assume that $C \in \on{GL}_n^{\mathrm{reg}}$.
Then the (commutative) algebra ${\bf{B}}(C)$ is a polynomial algebra in the coefficients of $\boldsymbol{\tau}_k(u,C)$.
Under the identification $\on{gr}{\bf{Y}}^{\mathrm{rtt}}(\mathfrak{gl}_n) \simeq \CO(\on{Mat}_n((z^{-1})))$ of Proposition \ref{iso_poiss_gr_yang_gl_univ}, we have $\on{gr}{\bf{B}}(C)=\ol{\bf{B}}(C)$.
\eprop
\prf
It follows from Proposition 
\ref{iso_poiss_gr_yang_gl_univ} that $\on{gr}\boldsymbol{\tau}_k(u,C)={\boldsymbol{\sigma}}_k(u,C)$, and thus the claim follows from Proposition \ref{size_of_univ_class_Bethe_gl}. 
\epr






\section{Antidominantly shifted Yangians for $\mathfrak{gl}_n$: standard and RTT realizations
}\label{shifted_yangians_section} 
We recall the RTT realization of $Y_\mu(\mathfrak{g})$  for classical $\mathfrak{g}$ that is given in \cite{fpt}, \cite{ft}.
As in Section \ref{univ_bethe_yang_sect!} we assume that 
$\mathfrak{g}=\mathfrak{gl}_n$.
Let $V=\BC^n$ be the tautological representation of $\mathfrak{g}$.  

Let us recall the shifted Yangian $Y_{\mu}(\mathfrak{gl}_n)$. Recall that $\epsilon_1,\ldots,\epsilon_n$ is the standard basis of $\BC^n$ and  $\epsilon_1^{\vee},\ldots,\epsilon_n^{\vee} \in (\BC^n)^*$ is the dual basis.
We consider the lattice $\La^{\vee}=\oplus_{j=1}^n \BZ \epsilon_j^{\vee}$ and the dual lattice $\La=\oplus_{j=1}^n \BZ \epsilon_j$. For $\mu \in \La$ we define 
\begin{equation*}
d_j:=\epsilon_j^{\vee}(\mu).
\end{equation*}


\begin{defeni}{}\label{BFN_def_shifted_yangian}
The shifted Yangian $Y_\mu(\mathfrak{gl}_n)$ is the associative $\BC$-algebra generated by $\{E_i^{(r)}, F_i^{(r)}\}^{r \geqslant 1}_{1 \leqslant i \leqslant n-1} \cup \{D_i^{(s_i)},\widetilde{D}_i^{(\widetilde{s}_i)}\}_{1 \leqslant i \leqslant n}^{s_i \geqslant -d_i, \widetilde{s}_i \geqslant d_i}$ with the following defining relations:
\begin{equation}
D_i^{(-d_i)}=1,\, \sum_{t=d_i}^{r+d_i}\widetilde{D}_i^{(t)}D_i^{(r-t)}=\delta_{r,0},\, [D_i^{(r)},D_j^{(s)}]=0,    
\end{equation}
\begin{equation}
[F_i^{(r)},E_j^{(s)}]=\delta_{i,j}\sum_{t=-d_i}^{r+s-1-d_{i+1}} D_i^{(t)} \widetilde{D}_{i+1}^{(r+s-t-1)},    
\end{equation}
\begin{equation}
[\widetilde{D}_i^{(r)},F_j^{(s)}]=(\delta_{i,j+1}-\delta_{i,j})\sum_{t=d_i}^{r-1}\widetilde{D}_i^{(t)}F_{j}^{(r+s-t-1)},
\end{equation}
\begin{equation}
[\widetilde{D}_i^{(r)},E_j^{(s)}]=(\delta_{i,j}-\delta_{i,j+1})\sum_{t=d_i}^{r-1}E_j^{(r+s-t-1)}\widetilde{D}_i^{(t)},    
\end{equation}
\begin{equation}
[F_i^{(r)},F_i^{(s)}]=\sum_{t=1}^{r-1}F_i^{(t)}F_i^{(r+s-t-1)}-\sum_{t=1}^{s-1}F_i^{(t)}F_i^{(r+s-t-1)},
\end{equation}
\begin{equation}
[E_i^{(r)},E_i^{(s)}]=\sum_{t=1}^{s-1}E_i^{(r+s-t-1)}E_i^{(t)}-\sum_{t=1}^{r-1}E_i^{(r+s-t-1)}E_i^{(t)},    
\end{equation}
\begin{equation}
[F_i^{(r+1)},F_{i+1}^{(s)}]-[F_i^{(r)},F_{i+1}^{(s+1)}]=-F_i^{(r)}F_{i+1}^{(s)},
\end{equation}
\begin{equation}
[E_i^{(r+1)},E_{i+1}^{(s)}]-[E_i^{(r)},E_{i+1}^{(s+1)}]=E_{i+1}^{(s)}E_i^{(r)},    
\end{equation}
\begin{equation}
[E_i^{(r)},E_j^{(s)}]=0~\text{if}~|i-j|>1, 
\end{equation}
\begin{equation}
[F_i^{(r)},F_j^{(s)}]=0~\text{if}~|i-j|>1,    
\end{equation}
\begin{equation}
[E_i^{(r)},[E_i^{(s)},E_j^{(t)}]]+[E_i^{(s)},[E_i^{(r)},E_j^{(t)}]]=0~\text{if}~|i-j|=1,    
\end{equation}
\begin{equation}
[F_i^{(r)},[F_i^{(s)},F_j^{(t)}]]+[F_i^{(s)},[F_i^{(r)},F_j^{(t)}]]=0~\text{if}~|i-j|=1.    
\end{equation}
We denote by $H \subset Y_\mu(\mathfrak{gl}_n)$ the (commutative) subalgebra generated by $\{D_i^{(s_i)},\widetilde{D}_i^{(\widetilde{s}_i)}\}_{1 \leqslant i \leqslant n}^{s_i \geqslant -d_i, \widetilde{s}_i \geqslant d_i}$ and call $H$ the {\emph{Cartan subalgebra}} of $Y_\mu(\mathfrak{gl}_n)$.
\end{defeni}

\rem{}
{\emph{Note that if we define 
\begin{equation*}
D_i(z)=\sum_{r \geqslant {{-}}d_i}D_i^{(r)}z^{-r},\, \widetilde{D}_i(z)=\sum_{r \geqslant d_i}\widetilde{D}_i^{(r)}z^{-r}    
\end{equation*}
then we have $D_i(z)\widetilde{D}_i(z)=1$.}} 
\erem

\begin{warning}{}
{\emph{The above definition of $Y_\mu(\mathfrak{gl}_n)$  slightly differs from the one given in \cite[Section 2.1]{fpt}.  If $Y^{\mathrm{FPT}}_\mu(\mathfrak{gl}_n)$ is the algebra defined in \cite[Section 2.1]{fpt} and $E_i^{\mathrm{FPT},(r)}$, $F_i^{\mathrm{FPT},(r)}$, $D_i^{\mathrm{FPT},(s_i)}$, $\widetilde{D}_i^{\mathrm{FPT},(\widetilde{s}_i)}$ are the generators of this algebra then the identification $Y_\mu(\mathfrak{gl}_n) \iso Y^{\mathrm{FPT}}_{\mu}(\mathfrak{gl}_n)$ is given by $E_{i}^{(r)} \mapsto -F_i^{\mathrm{FPT},(r)}$, $F_i^{(r)} \mapsto -E_i^{\mathrm{FPT},(r)}$, $D_i^{(s_i)} \mapsto -\widetilde{D}^{\mathrm{FPT},(s_i)}$, $\widetilde{D}^{(\widetilde{s}_i)} \mapsto D^{\mathrm{FPT},(\widetilde{s}_i)}_i$, for $r \geqslant 1$, $s_i \geqslant -d_i$, $\widetilde{s}_i \geqslant d_i$.
}}
\end{warning}

Similarly to \cite[Section 5.4]{fkprw} we consider the filtration on $Y_\mu(\mathfrak{gl}_n)$ such that 
\begin{equation}\label{def_filtr_gen}
\on{deg}E_{i}^{(r)}=r,\, \on{deg}F_i^{(r)}=r+d_i-d_{i+1},\, \on{deg}D_i^{(s_i)}=s_i+d_i,\, \on{deg}\widetilde{D}_i^{(\tilde{s}_i)}=\tilde{s}_i-d_i.
\end{equation}
\begin{warning}{}\label{warn_filtr_pbw}
The assignment (\ref{def_filtr_gen}) is {\emph{insufficient}} to define our filtration. One should consider a  ``PBW'' basis in $Y_\mu(\mathfrak{gl}_n)$ and define the degree of each of the basis elements (as \cite{fkprw} do in $(5.1)$). The relevant PBW basis is described in Proposition \ref{PBW_for_Y_mu} below, and the filtration is defined in (\ref{def_filtr_basis}).   
\end{warning}


Recall the Poisson subvariety $\CW_\mu \subset \on{GL}_n((z^{-1})) 
 \subset \on{Mat}_n((z^{-1}))$ (see (\ref{W_mu_for_antidom_gl}), (\ref{gen_def_W_mu})) that can be presented as 
\begin{equation*}
\CW_\mu = U[[z^{-1}]]_1 z^\mu T[[z^{-1}]]_1 U_-[[z^{-1}]]_1. 
\end{equation*}
Recall also that we have a $\BC^\times$-action on $\CW_\mu$ given by $t \cdot g(z)=g(tz)t^{-\mu}$.
This action induces the grading on $\CO(\CW_\mu)$.

Let $\ol{e}_{ij}^{(r)} \in \CO(\CW_\mu)$ be the function that sends $g=u \cdot z^\mu \cdot t \cdot u_-$ to the $z^{-r}$-coefficient of $(i,j)$th matrix coefficient of $u$, let $\ol{f}_{ji}^{(r)} \in \CO(\CW_\mu)$ be the function that sends $g=u \cdot z^\mu \cdot t \cdot u_-$ to the $z^{-r}$-coefficient of $(j,i)$th matrix coefficient of $u_-$, and let $\ol{g}_j^{(r)}$ be the function that sends $g=u \cdot z^\mu \cdot t \cdot u_-$ to the $z^{-r}$-coefficient of $(j,j)$th matrix coefficient of $z^\mu \cdot t$. Let $\ol{F}(u), \ol{G}(u), \ol{E}(u) \in \CO(\CW_\mu)((u^{-1})) \otimes \on{End}(\BC^n)$ be the corresponding generating functions. We note that they are lower-triangular, diagonal, and upper-triangular matrices.

The following proposition is similar to  \cite[Theorem 5.15]{fkprw} together with \cite[Theorem A.11]{kpw}.
\prop{}\label{gr_of_shift_yang_expl}
There exists the isomorphism of Poisson graded  algebras $\on{gr}Y_\mu(\mathfrak{gl}_n) \simeq \CO(\CW_\mu)$, which identifies $\on{gr}E_i^{(r)}$ with the function $\ol{e}_{i\,i+1}^{(r)}$, $\on{gr}F_{i}^{(r)}$ with the function $\ol{f}_{i+1\,i}^{(r)}$, and $\on{gr}D_i^{(r)}$ with the function $\ol{g}_i^{(r)}$.
\eprop

We shall
slightly modify
the notations of \cite{fpt}.
Recall  the rational $R$-matrix $R(u)$:
\begin{equation*}
R(u)=-u - P,\, P=\sum_{i,j}E_{ij} \otimes E_{ji},   
\end{equation*}
where $E_{ij} \in \mathfrak{gl}_n$ are the matrix units. 


\begin{defeni}{}\label{defe_RTT_Yang!}
The (antidominantly) shifted RTT Yangian $Y^{\mathrm{rtt}}_\mu(\mathfrak{gl}_n)$ is the associative $\BC$-algebra generated by $\{t_{ij}^{(r)}\}^{r \in \BZ}_{1 \leqslant i,j \leqslant n}$ subject to the following two families of relations:
\begin{equation*}
R(u_1-u_2)T_1(u_1)T_2(u_2)=T_2(u_2)T_1(u_1)R(u_1-u_2),    
\end{equation*}
where $T(u) \in Y^{\mathrm{rtt}}_\mu(\mathfrak{gl}_n)[[u,u^{-1}]] \otimes \on{End}(\BC^n)$ is defined via 
\begin{equation*}
T(u)=(t_{ij}(u))_{ij}=\sum_{ij}t_{ij}(u) \otimes E_{ij}~\text{with}~t_{ij}(u)=\sum_{r \in \BZ} t_{ij}^{(r)}u^{-r}.    
\end{equation*}
The second family of relations encodes the fact that $T(u)$ admits the Gauss decomposition:
\begin{equation}\label{Gauss_decomp}
T(u)=E(u) \cdot G(u) \cdot F(u),    
\end{equation} 
where $E(u), G(u), F(u) \in Y^{\mathrm{rtt}}_\mu(\mathfrak{gl}_n)((u^{-1})) \otimes \on{End}(\BC^n)$ are of the form
\begin{equation*}
E(u)=\sum_i E_{ii}+\sum_{i < j} e_{ij}(u) \otimes E_{ij},\, G(u)=\sum_{i} g_i(u) \otimes E_{ii},\,  F(u)=\sum_i E_{ii}    + \sum_{i<j}f_{ji}(u) \otimes E_{ji}
\end{equation*}
with the matrix coefficients having the following expansions in $u$:
\begin{equation*}
e_{ij}(u)=\sum_{r \geqslant 1} e_{ij}^{(r)}u^{-r},\, f_{ji}(u)=\sum_{r \geqslant 1}f_{ji}^{(r)}u^{-r},\, g_i(u)=u^{d_i}+\sum_{r \geqslant 1-d_i}g_i^{(r)}u^{-r}, 
\end{equation*}
where $d_i=\epsilon_i^{\vee}(\mu)$ as before.
\end{defeni}






\begin{warning}{}
{\emph{In \cite{fpt} authors consider the $R$-matrix $R^{\mathrm{FPT}}(u)=-R(u)$. Let $Y^{\mathrm{FPT},\mathrm{rtt}}_\mu(\mathfrak{gl}_n)$ be the RTT Yangian as in \cite[Section 2.3]{fpt}. If $t_{ij}^{\mathrm{FPT},(r)}$, $e_{ij}^{\mathrm{FPT},(r)}$, $f_{ji}^{\mathrm{FPT},(r)}$, $g_i^{\mathrm{FPT},(r)}$ are  the elements as in \cite[Section 2.3]{fpt} and $T^{\mathrm{FPT}}(u)$, $E^{\mathrm{FPT}}(u)$, $F^{\mathrm{FPT}}(u)$, $G^{\mathrm{FPT}}(u)$ are the corresponding matrices  then we have the isomorphism of algebras $Y^{\mathrm{rtt}}_\mu(\mathfrak{gl}_n) \iso Y^{\mathrm{FPT},\mathrm{rtt}}_\mu(\mathfrak{gl}_n)$ given by 
\begin{equation*}
T(u) \mapsto (T^{\mathrm{FPT}}(u)^{-1})^t,
\end{equation*}
\begin{equation*}
E(u) \mapsto (F^{\mathrm{FPT}}(u)^{-1})^t,\, F(u) \mapsto (E^{\mathrm{FPT}}(u)^{-1})^t,\, G(u) \mapsto G^{\mathrm{FPT}}(u)^{-1}
\end{equation*}
(similarly to \cite[Propositions 1.11, 1.12 and Corollary 1.13]{mno} one can see that this assignment indeed extends to the isomorphism of algebras).
Note that this isomorphism sends $e^{(r)}_{i,i+1}$ to $-f^{\mathrm{FPT},(r)}_{i+1,i}$, $f^{(r)}_{i+1,i}$ to $-e^{\mathrm{FPT},(r)}_{i,i+1}$.
We will see below (Lemma \ref{gr_t_i_j!}) that $Y_{\mu}(\mathfrak{gl}_n)$ admits a filtration such that we have an isomorphism of Poisson algebras $\on{gr}Y_{\mu}(\mathfrak{gl}_n) \simeq \CO(\CW_\mu)=\CO(U[[z^{-1}]]_1 z^\mu T[[z^{-1}]]_1 U_-[[z^{-1}]])$ given by $\on{gr}E(u) \mapsto \ol{E}(u)$, $\on{gr}G(u) \mapsto \ol{G}(u)$, $\on{gr}F(u) \mapsto \ol{F}(u)$. This is the main motivation for us to slightly change the RTT definition of the shifted Yangian given in \cite[Section 2.3]{fpt}. Note that $\on{gr}Y_\mu^{\mathrm{FPT},\mathrm{rtt}}(\mathfrak{gl}_n)$ is naturally identified with $\CO(U_-[[z^{-1}]]_1 z^{-\mu} T[[z^{-1}]]_1 U[[z^{-1}]])$.
}}
\end{warning}


The following lemma is clear from the Gauss decomposition (\ref{Gauss_decomp}):
\lem{}
For every $i,j=1,\ldots,n$,
we have $t_{ij}^{(r)}=0$ for $r \ll 0$.
\elem

As a corollary, we obtain:
\cor{}\label{surj_univ_yang_to_shifted}
We have a surjective homomorphism of algebras ${\bf{Y}}^{\mathrm{rtt}}(\mathfrak{gl}_n) \twoheadrightarrow Y^{\mathrm{rtt}}_\mu(\mathfrak{gl}_n)$ given by ${\bf{T}}(u) \mapsto T(u)$.
\ecor

The following Lemma holds.
\lem{}\label{realiz_e_from_smal}$(a)$  For $i+1<j \leqslant n$, 
we have the following identities in $Y_\mu^{\mathrm{rtt}}(\mathfrak{gl}_n)$: 
\begin{equation*}
e_{ij}^{(r)}=[e_{i+1,j}^{(r)},e^{(1)}_{i,i+1}],\, f_{ji}^{(r)}=[f_{i+1,i}^{(1)},f_{j,i+1}^{(r)}]. 
\end{equation*}

$(b)$ For $i+1<j \leqslant n$, 
we have the following identities in $\CO(\CW_\mu)$:
\begin{equation*}
\ol{e}_{ij}^{(r)}=\{\ol{e}_{i+1,j}^{(r)},\ol{e}^{(1)}_{i,i+1}\},\, \ol{f}_{ji}^{(r)}=\{\ol{f}_{i+1,i}^{(1)},\ol{f}_{j,i+1}^{(r)}\}.
\end{equation*}
\elem
\prf
Part $(a)$ is similar to \cite[Lemma 2.46]{fpt} (see {\cite[Equation (5.5)]{bk0}}), the proof of part $(b)$ is analogous.
\epr

The following theorem is proven in \cite[Theorem 2.54]{fpt}. 
\begin{thm}{}\label{rtt_iso_Drinf_thm}
For any antidominant $\mu \in \La$ there exists the unique isomorphism of algebras $\Theta\colon Y_\mu^{\mathrm{rtt}}(\mathfrak{gl}_n) \iso Y_\mu(\mathfrak{gl}_n)$ given by
\begin{equation*}
   e_{i,i+1}(u) \mapsto E_i(u) ,\, f_{i+1,i}(u) \mapsto F_i(u) ,\, g_j(u) \mapsto D_j(u).  
\end{equation*}
\end{thm}



From now on we identify the algebras $Y^{\mathrm{rtt}}_{\mu}(\mathfrak{gl}_n)=Y_\mu(\mathfrak{gl}_n)$ via $\Theta$.
The following proposition follows from \cite[Corollary 2.24]{fpt}, \cite[Theorem 3.14]{fkprw}, and Theorem \ref{rtt_iso_Drinf_thm}.

\prop{}\label{PBW_for_Y_mu} (PBW theorem for $Y_\mu(\mathfrak{gl}_n)$)
The set of ordered monomials in the variables 
\begin{equation*}
\{e_{ij}^{(r)},\, i<j, r>0\} \cup \{g_i^{(s_i)},\,s_i>-d_i\}  \cup \{f_{ji}^{(r)},\, i<j, r>0\} 
\end{equation*}
forms a basis of $Y_\mu(\mathfrak{gl}_n)$ over $\BC$.  The filtration (\ref{def_filtr_gen}) above can be defined as follows (cf. Warning \ref{warn_filtr_pbw}): 
\begin{equation}\label{def_filtr_basis}
\on{deg}e_{ij}^{(r)}=r,\,\on{deg}g_i^{(s_i)}=s_i+d_i,\, \on{deg}f_{ji}^{(r)}=r+d_i-d_j.    
\end{equation}
\eprop

Let us recall the isomorphism of Poisson algebras $\on{gr}Y_\mu(\mathfrak{gl}_n) \simeq \CO(\CW_\mu)$ constructed in Proposition \ref{gr_of_shift_yang_expl}.
\lem{}\label{ident_higher_e_f}
The isomorphism $\on{gr}Y_\mu(\mathfrak{gl}_n) \simeq \CO(\CW_\mu)$  identifies $\on{gr}e_{ij}^{(r)}$ with the function $\ol{e}_{ij}^{(r)}$, $\on{gr}f_{ji}^{(r)}$ with the function $\ol{f}_{ji}^{(r)}$ and $\on{gr}D_i^{(r)}$ with the function $\ol{g}_i^{(r)}$.
\elem
\prf
Follows from Proposition \ref{gr_of_shift_yang_expl}, together with Lemma \ref{realiz_e_from_smal}.
\epr

\lem{}\label{gr_t_i_j!}
The element $t_{ij}^{(r)} \in Y^{\mathrm{rtt}}_\mu(\mathfrak{gl}_n)$ has degree $r+d_j$ with respect to the filtration above. In particular, the homomorphism ${\bf{Y}}^{\mathrm{rtt}}(\mathfrak{gl}_n) \ra Y_\mu(\mathfrak{gl}_n)$ becomes a homomorphism of filtered algebras if we endow ${\bf{Y}}^{\mathrm{rtt}}(\mathfrak{gl}_n)$ with the ``$\mu$-twisted'' filtration (\ref{filtr_twist_univ}). We have 
$\on{gr}t_{ij}^{(r)}=\Delta_{ij}^{(r)} \in \CO(\CW_\mu)$. 
\elem
\prf 
The equality $T(u)=E(u) \cdot G(u) \cdot F(u)$ can be rewritten as 
\begin{equation*}
t_{ij}(u)=\sum_{k \geqslant \on{max}(i,j)} e_{ik}(u)g_{k}(u)f_{kj}(u),\, i,j=1,\ldots,n,    
\end{equation*}
where we assume that $e_{kk}(u)=f_{kk}(u)=1$ for every $k=1,\ldots,n$. 
So 
\begin{equation*}
t_{ij}^{(r)}=\sum_{k \geqslant \on{max}(i,j),\, r_1+r_2+r_3=r} e_{ik}^{(r_1)}g_k^{(r_2)}f_{kj}^{(r_3)}.  
\end{equation*}
Note now that $\on{deg}(e_{ik}^{(r_1)}g_k^{(r_2)}f_{kj}^{(r_3)})=r_1+r_2+d_k+r_3+d_j-d_k=r+d_j$. Since the elements $e_{ik}^{(r_1)}g_k^{(r_2)}f_{kj}^{(r_3)}$ form a subset of the PBW-basis for $Y_\mu(\mathfrak{gl}_n)$ (see Proposition \ref{PBW_for_Y_mu} above)  we conclude that $\on{deg}t_{ij}^{(r)}=r+d_j$.

It remains to show that $\on{gr}t_{ij}^{(r)}=\Delta_{ij}^{(r)}$. 
Recall that $\ol{T}(u)=\ol{E}(u) \cdot \ol{G}(u) \cdot \ol{F}(u)$, the matrix coefficients of $\ol{T}(u)$ are $\Delta_{ij}(u)=\sum_{r \in \BZ} \Delta_{ij}^{(r)}u^{-r}$.  It follows from Lemma \ref{ident_higher_e_f} that 
\begin{equation*}
\on{gr}E(u)=\ol{E}(u),\, \on{gr}G(u)=\ol{G}(u),\, \on{gr}F(u)=\ol{F}(u).    
\end{equation*} Using the decomposition $T(u)=E(u) \cdot G(u) \cdot F(u)$ we conclude that $\on{gr}t_{ij}^{(r)}=\Delta_{ij}^{(r)}.$

\epr

We obtain the following corollary.
\cor{}\label{surj_quant_surj}
The natural surjection ${\bf{Y}}^{\mathrm{rtt}}(\mathfrak{gl}_n) \twoheadrightarrow Y_\mu(\mathfrak{gl}_n)$ quantizes the surjection $\CO(\on{Mat}_n((z^{-1}))) \twoheadrightarrow \CO(\CW_\mu)$ induced by the embedding $\CW_\mu \subset \on{Mat}_n((z^{-1}))$.   
\ecor
\prf
It follows from Lemma \ref{gr_t_i_j!} that the surjection ${\bf{Y}}^{\mathrm{rtt}}(\mathfrak{gl}_n) \twoheadrightarrow Y_\mu(\mathfrak{gl}_n)$ is indeed a homomorphism of filtered algebras (recall that we consider the ``$\mu$-twisted'' filtration on ${\bf{Y}}^{\mathrm{rtt}}(\mathfrak{gl}_n)$). It follows from Proposition \ref{iso_poiss_gr_yang_gl_univ_twisted} that the identification $\on{gr}_\mu{\bf{Y}}^{\mathrm{rtt}}(\mathfrak{gl}_n) \iso \CO(\on{Mat}_n((z^{-1})))$ sends $\on{gr}_\mu t_{ij}^{(r)}$ to $\Delta_{ij}^{(r)}$. Recall also that by Lemma \ref{gr_t_i_j!} the associated graded of $t_{ij}^{(r)} \in Y_\mu(\mathfrak{gl}_n)$ is equal to $\Delta_{ij}^{(r)} \in \CO(\CW_\mu)$.  We conclude that the associated graded of the surjection ${\bf{Y}}^{\mathrm{rtt}}(\mathfrak{gl}_n) \twoheadrightarrow Y_\mu(\mathfrak{gl}_n)$ induces the homomorphism  $\CO(\on{Mat}_n((z^{-1}))) \ra \CO(\CW_\mu)$ that sends $\Delta_{ij}^{(r)}$ to its restriction to $\CW_\mu$. It follows that this map is induced by the natural (closed) embedding $\CW_\mu \subset \on{Mat}_n((z^{-1}))$.
\epr

\section{Bethe subalgebras in $Y_\mu(\mathfrak{gl}_n)$}\label{bethe_subalg_in_shifted_section!!!}

\subsection{Commutative subalgebras $B_\mu(C)$}{} Recall that for $k=1,\ldots,n$ we denote by $A_k \in \on{End}(\BC^n)^{\otimes k}$ the antisymmetrization  map normalized so that $A_k^2=A_k$. We set
\begin{equation*}
\tau_{\mu,k}(u,C):=\on{tr}_{(\BC^n)^{\otimes k}} A_k C_1 \ldots C_k T_1(u) \ldots T_k(u-k+1) \in Y^{\mathrm{rtt}}_\mu(\mathfrak{gl}_n)((u^{-1})).   
\end{equation*}

As proposed in \cite{fpt}, we make:
\begin{defeni}{}\label{defe_bethe_shifted}
We denote by $B_\mu(C) \subset Y_\mu^{\mathrm{rtt}}(\mathfrak{gl}_n)$ the subalgebra generated by the coefficients of $\tau_{\mu,k}(u,C)$, $k=1,\ldots,n$.
\end{defeni}

Recall now that in Section \ref{univ_bethe_in_yang} we have defined ``universal'' Bethe subalgebra ${\bf{B}}(C) \subset {\bf{Y}}^{\mathrm{rtt}}(\mathfrak{gl}_n)$ and proved that it is commutative. Recall also that we have the natural surjection (see Corollary \ref{surj_univ_yang_to_shifted}) \begin{equation*}
{\bf{Y}}^{\mathrm{rtt}}(\mathfrak{gl}_n) \twoheadrightarrow Y_\mu(\mathfrak{gl}_n)    
\end{equation*}
that sends ${\bf{B}}(C)$ onto $B_\mu(C)$,  as follows from the definitions.
In particular, we recover the observation of \cite{fpt}:
\prop{}\label{shifted_bethe_are_comm}
The algebra $B_\mu(C)$ is commutative.
\eprop

\subsection{Size of $B_\mu(C)$ and the associated graded $\on{gr}B_\mu(C)$}{}\label{size_shift_bethe}  Let $\tau_{\mu,k}(C)^{(r)}$ be the coefficient of $u^{-r}$ in $\tau_{\mu,k}(u,C)$. 
Recall that we have a filtration on the algebra $Y_\mu(\mathfrak{gl}_n)=Y^{\mathrm{rtt}}_{\mu}(\mathfrak{gl}_n)$ such that the associated graded algebra is isomorphic to $\CO(\CW_\mu)$. Recall also the functions $\sigma_{\mu,k}(u,C)$, $k=1,\ldots,n$ defined  in Section \ref{reform_for_gl}. 
For $k=1,\ldots,n$ recall that $\omega_k^*=-\epsilon_n^\vee-\ldots-\epsilon_{n-k+1}^\vee$. 


\lem{}\label{zero_tau_for_small}
The element $\tau_{\mu,k}(C)^{(r)}$ is equal to zero for $r < \langle \omega_k^*,\mu \rangle$ and $\tau_{\mu,k}(C)^{(\langle \omega_k^*,\mu \rangle)}$ is a positive integer.
\elem
\prf
The proof is the same as the one of Lemma \ref{sigma_zero_neg_degree}: one should use the decomposition $T(u)=E(u) \cdot G(u) \cdot F(u)$ and the fact that $\tau_{\mu,k}(C,u)$ is the linear combination of the elements of the form $t_{a_1b_1}(u)\ldots t_{a_kb_k}(u-k+1)$ where $a_i \neq a_j$, $b_i \neq b_j$ for $1 \leqslant i < j \leqslant n$ and $\{a_1,a_2,\ldots,a_k\}=\{b_1,b_2,\ldots,b_k\}$ (use that $A_k$ is the projector onto $\La^k(\BC^n) \subset (\BC^n)^{\otimes k}$). 
\epr

Recall that $L=Z_{G}(\mu)$. 

\th{}\label{thm_size_shifted}
For $C \in L^{\mathrm{reg}}$, $k=1,\ldots,n$, and $r>\langle \omega_k^*,\mu \rangle$ the element $\tau_{\mu,k}(C,u)^{(r)}$ has degree $r-\langle \omega_k^*,\mu \rangle$, and we have 
\begin{equation*}
\on{gr}\tau_{\mu,k}(C)^{(r)}=\on{gr}\sigma_{\mu,k}(C)^{(r)}.
\end{equation*} 
\eth
\prf
Recall that 
\begin{equation*}
\tau_{\mu,k}(C,u)= \on{tr}_{(\BC^n)^{\otimes k}} A_k C_1\ldots C_k T_1(u) \ldots T_k(u-k+1).   
\end{equation*}
Note also that 
\begin{equation}
\on{tr}_{(\BC^n)^{\otimes k}}A_k C_1 \ldots C_k \ol{T}_1(u) \ldots \ol{T}_k(u)=\sigma_{\mu,k}(C,u)    
\end{equation}
($\sigma_{\mu,k}(C,u)$ is defined in Definition \ref{class_bethe_in_shifted_gl}).
It is  easy to conclude from Lemma \ref{gr_t_i_j!} that for every $r \in \BZ$ the element $\tau_{\mu,k}(C)^{(r)}$ has degree at most $r-\langle \omega_k^*,\mu\rangle$. Recall now the element $\sigma_{\mu,k}(C)^{(r)}$ (see Definition \ref{class_bethe_in_shifted_gl}). It follows from Theorem \ref{size_reg_bethe_gl} that the degree of this element is {\em{equal}} to $r-\langle \omega_k^*,\mu\rangle$. This means that $\on{gr}\sigma_{\mu,k}(C)^{(r)} \in \CO(\CW_\mu)$ is the nonzero element of degree $r-\langle \omega_k^*,\mu\rangle$.
We now conclude from Lemma \ref{gr_t_i_j!} and from the definition of $\tau_{\mu,k}(C)^{(r)}$ that  $\tau_{\mu,k}(C)^{(r)}$ is the nonzero element of degree $r-\langle \omega_k^*,\mu\rangle$ so it is clear that $\on{gr} \tau_{\mu,k}(C)^{(r)}=\on{gr}\sigma_{\mu,k}(C)^{(r)}$ (equality of the elements of $\on{gr}\CO(\CW_\mu) \simeq \CO(\CW_\mu)$).

\epr


We are now ready to prove Theorem \ref{thm_C}. 

\cor{}\label{cor_size_Bethe_type_A}
The composition 
\begin{equation*}
\on{gr}B_\mu(C) \hookrightarrow \CO(\CW_\mu) \twoheadrightarrow \CO(L[[z^{-1}]]_1)    
\end{equation*}
induces the identification 
\begin{equation*}
\on{gr}B_\mu(C) \iso \ol{B}_L(C)    
\end{equation*}
which sends $\on{gr}\tau_{\mu,k}(C)^{(r)}$ to $\sigma_{L,k}(C)^{(r-\langle \omega_k^*,\mu \rangle)}$.
The algebra $B_\mu(C)$ is a polynomial algebra in the elements $\{\tau_{\mu,k}(C)^{(r)},\, r > \langle \omega_k^*,\mu \rangle\}$ and the Poincar\'e series of $B_\mu(C)$ coincides with the Poincar\'e series of the Cartan subalgebra $H \subset Y_\mu(\mathfrak{gl}_n)$.
\ecor
\prf
Follows from Corollary \ref{cor_size_class_Bethe_shifted_gl}, Lemma~\ref{zero_tau_for_small}, and Theorem \ref{thm_size_shifted}.
\epr

\appendix
\section{Comparing Bethe subalgebras in $\CO(\CW_\mu)$ with pullbacks of Bethe subalgebras in $\CO(\CW_0)$}\label{app_1} 
In this section, we assume that $\mathfrak{g}=\mathfrak{sl}_2$. Recall that $T \subset \on{SL}_2$ is the subgroup of diagonal matrices and $B \subset \on{SL}_2$ is the subgroup of upper triangular matrices. We identify $\on{Hom}(\BC^\times,T)$ with $\BZ$ and fix $\mu=-n$ for some $n \in \BZ_{\geqslant 0}$.
By \cite[Section 2(xii)]{bfn} we have:
\begin{multline*}
\CW_{-n}=\\
=\Big\{
\begin{pmatrix}
a & b \\
c & d
\end{pmatrix}\,|\, a,b,c,d \in \BC((z^{-1})),\, d=z^n+\ldots,\, n > \text{top}~z\text{-degree of the series}~b, c,\, ad-bc=1\Big\},    
\end{multline*}
in particular $\CW_0=\on{SL}_2[[z^{-1}]]_1$.

It is easy to see that the shift morphism $\iota_{0,0,-n}$ (from \cite[Section 5.9]{fkprw}) sends 
\begin{equation*}
g=
\begin{pmatrix}
xz^{-n}+x^{-1}QP z^n & x^{-1}Qz^n \\
x^{-1}P z^n & x^{-1}z^n
\end{pmatrix}=
\begin{pmatrix}
1 & Q \\
0 & 1
\end{pmatrix}\begin{pmatrix}
x & 0 \\
0 & x^{-1}
\end{pmatrix}
\begin{pmatrix}
z^{-n} & 0 \\
0 & z^n
\end{pmatrix}\begin{pmatrix}
1 & 0\\
P & 1
\end{pmatrix}
\end{equation*}
to
\begin{equation*}
\begin{pmatrix}
x+x^{-1}QP' & x^{-1}Q \\
x^{-1}P' & x^{-1}
\end{pmatrix}=
\begin{pmatrix}
1 & Q \\
0 & 1
\end{pmatrix}\begin{pmatrix}
x & 0 \\
0 & x^{-1}
\end{pmatrix}
\begin{pmatrix}
1 & 0\\
P' & 1
\end{pmatrix},    
\end{equation*}
where $P'=\sum_{k=1}^{\infty}z^{-k}P_{-k-2n}$.

Recall now that to $C, D \in \on{SL}_2$ one can associate subalgebras $\ol{B}_{-n}(C) \subset \CO(\CW_{-n})$, $\ol{B}(D) \subset \CO(\on{SL}_2[[z^{-1}]]_1)$ and our goal is to compare $\ol{B}_{-n}(C)$ with $\iota_{0,0,-n}^*(\ol{B}(D))$. 

\prop{}\label{alg_are_diff_prf}
For $n>0$ and $C \in T$, and any $D \in \on{SL}_2$ the algebras $\ol{B}_{-n}(C)$, $\iota_{0,0,-n}^*(\ol{B}(D))$ are distinct. \eprop
\prf
Let $C=\on{diag}(h,h^{-1})$ for some $h \in \BC^\times$.
Recall that 
$
\CW_{-n}=B[[z^{-1}]]_1z^{\mu}U_-[[z^{-1}]]_1    
$
and consider 
$Z:=\Big\{\begin{pmatrix}
1 & 
z^{-1} \\
0 & 1
\end{pmatrix}
\cdot z^\mu \cdot \begin{pmatrix}
1 & 0 \\
\kappa z^{-1} & 1
\end{pmatrix}\,|\, \kappa \in \BC\Big\} \subset \CW_{-n}$. 

Recall now that the algebra $\ol{B}_{-n}(C)$, $C=\on{diag}(h,h^{-1})$ is generated by the coefficients of the function $g \mapsto \on{tr}(Cg)$. Being restricted to $Z$ this function is given by 
\begin{equation*}
\kappa \mapsto hz^{-n}+h
\kappa z^{n-2}+h^{-1}z^n. \end{equation*}
So the image of $\ol{B}_{-n}(C)$ in $\CO(Z)$ will be $\BC[\kappa]$.

Note now that
\begin{equation*}
\iota_{0,0,-n}\Big(\begin{pmatrix}
1 & 
z^{-1} \\
0 & 1
\end{pmatrix} \begin{pmatrix}
z^{-n} & 0 \\
0 & z^n
\end{pmatrix}
\begin{pmatrix}
1 & 0 \\
\kappa z^{-1} & 1
\end{pmatrix}
\Big)=\begin{pmatrix}
1 & 
z^{-1} \\ 
0 & 1
\end{pmatrix}
\end{equation*}
so every function of the form $\iota_{0,0,-n}^*(f)$ ($f \in \CO(\CW_0)$) is constant on $Z$.
We conclude that the images of $\iota_{0,0,-n}^*(\ol{B}(D))$, $\ol{B}_{-n}(C)$ in $\CO(Z)$ are different algebras. 
\epr

\rem{}
{\em{Similar argument as the one in Proposition \ref{alg_are_diff_prf} applies for more general shift morphisms, where $\mu$ is split into two terms $\mu_1,\mu_2$ (see \cite[Proposition 3.8]{fkprw}). There exists the  ``Yangian'' (i.e. quantized) version of Proposition \ref{alg_are_diff_prf}, we omit the details.}}
\erem

\section{Generalization of a theorem of Steinberg}\label{app}
In this section, we formulate and prove a generalization of the classical theorem (due to Steinberg) that claims that the differentials at a regular element of characters of fundamental representations of semisimple simply connected groups are linearly independent. We generalize this result to the case of a Levi subgroup $L$ of a reductive group $G$ such that $[G,G]$ is simply connected. Our generalization is well-known to the experts, but we have decided to include the proof for completeness.  


\subsection{}
Let $L$ be {\em{any}} reductive group over complex numbers. 
Let $T \subset L$ be a maximal torus. Recall that an element $C \in L$ is called regular if $\on{dim}Z_L(C)=\on{dim}T$. Let $\on{Rep}L$ be the category of finite dimensional representations of $L$.
For $V \in \on{Rep}L$ let $\chi_{L,V} \colon L \ra \BC$ be its character.
Let $K(\on{Rep}L)$ be the Grothendieck group of the abelian category $\on{Rep}L$.
We have the ring homomorphism $K(\on{Rep}L) \ra \CO(L)^L$ given by $[V] \mapsto \chi_{L,V}$. The following proposition is standard. 
\prop{}\label{thre_descr_inv}
We have isomorphisms of algebras:
\begin{equation*}
K(\on{Rep}L) \otimes_{\BZ} \BC \iso \CO(L)^L \iso \CO(T)^W.
\end{equation*}
\eprop
\prf 
Standard, see for example \cite[Theorem 4]{ser}.
\epr

Let $\varphi\colon L \ra T/W$ be the morphism induced by the embedding 
\begin{equation*}
\CO(T)^W \simeq \CO(L)^L \subset \CO(L).
\end{equation*}

\prop{}\label{kernel_varphi}
The morphism $\varphi|_{L^{\mathrm{reg}}} \colon L^{\mathrm{reg}} \ra T/W$ is smooth. In particular, for every $C \in L^{\mathrm{reg}}$ we have an exact sequence 
\begin{equation*}
0 \ra T_C(L \cdot C) \ra T_CL \ra T_{\varphi(C)}(T/W) \ra 0.
\end{equation*}
\eprop
\prf
Let $\mathcal{D}(L):=[L,L]$ be the derived subgroup of $L$. Let $\widetilde{\mathcal{D}(L)}$ be the simply connected cover of $\mathcal{D}(L)$. Let $Z(L) \subset L$ be the center of $L$ and let $Z(L)^{\circ} \subset Z(L)$ be the connected component of $1 \in Z(L)$. Set $\tilde{L}:=\widetilde{\mathcal{D}(L)} \times Z(L)^{\circ}$.
Consider the natural central isogeny $\tilde{L} \ra L$ (compare with \cite[Theorem 3.2.2]{konrad}). Let $\tilde{T} \subset \tilde{L}$ be the preimage of $T \subset L$, $\tilde{T}$ is a maximal torus of $\tilde{L}$. It follows from 
\cite[Section 3.8]{st}  that the natural morphism $\tilde{L}^{\mathrm{reg}} \ra \tilde{T}/W$ is smooth. Set $T^{\mathrm{reg}}:=L^{\mathrm{reg}} \cap T$, $\tilde{T}^{\mathrm{reg}}=\tilde{L}^{\mathrm{reg}} \cap \tilde{T}$.
Recall now that we have the  cartesian diagram 
\begin{equation*}
\xymatrix{
\tilde{L}^{\mathrm{reg}} \ar[d] \ar[r]  & L^{\mathrm{reg}} \ar[d] \\ 
\tilde{T}^{\mathrm{reg}}/W \ar[r] & T^{\mathrm{reg}}/W
}
\end{equation*}
and the morphism $\tilde{T}^{\mathrm{reg}}/W \ra T^{\mathrm{reg}}/W$ is \'etale and surjective. We conclude that $L^{\mathrm{reg}} \ra T^{\mathrm{reg}}/W$ is smooth (use smooth descent w.r.t. flat, surjective morphism), so $L^{\mathrm{reg}} \ra T/W$ is smooth (morphism $T^{\mathrm{reg}}/W \rightarrow T/W$ is an open embedding, hence, smooth).
\epr


\prop{}\label{gen_lie_centr}
For $C \in L^{\mathrm{reg}}$ the restrictions of differentials \begin{equation*}
\{d_C(\chi_{L,V})\,|\, V \in \on{Rep}(L)\}
\end{equation*}
to $T_C Z_L(C)$ generate $T^*_C Z_L(C)$.
\eprop
\prf
We have an exact sequence 
\begin{equation*}
0 \ra T_C(L \cdot C) \ra T_C L \ra T_C (Z_L(C)) \ra 0,
\end{equation*}
where $L \cdot C$ is the $L$-orbit of $C \in L$ w.r.t. the adjoint action.
From Proposition \ref{kernel_varphi} we have an exact sequence 
\begin{equation*}
0 \ra T_C (L \cdot C) \ra T_C L \ra T_{\varphi(C)}(T/W) \ra 0.
\end{equation*}
We obtain the natural identification $T_C(Z_L(C)) \simeq T_{\varphi(C)}(T/W)$. It follows from the definitions that $T_{\varphi(C)}^*(T/W)$ is generated by $\{d_{\varphi(C)}f\,|\, f \in \CO(T)^W\}$. Now the claim follows from Proposition \ref{thre_descr_inv}. 
\epr

\subsection{} Assume now that $\mathcal{D}(L):=[L,L]$ is simply connected and set $\mathcal{D}(\mathfrak{l}):=[\mathfrak{l},\mathfrak{l}]$. Let $V_{L,1},\ldots, V_{L,\on{rk}\mathcal{D}(\mathfrak{l})}$ be irreducible (finite dimensional) modules over $L$ that restrict to fundamental representations of $\mathcal{D}(L)$. For $i=1,\ldots,\on{rk}\mathcal{D}(\mathfrak{l})$ let $\la_i \colon T \ra \BC^\times$ be the highest weight of $V_{L,i}$. Let $\BC_{\nu_1},\ldots,\BC_{\nu_{\on{rk}\mathfrak{l}-\on{rk}\mathcal{D}(\mathfrak{l})}}$ be one-dimensional representations of $L$ with highest weights $\nu_1,\ldots,\nu_{\on{rk}\mathfrak{l}-\on{rk}\mathcal{D}(\mathfrak{l})}$ such that the characters $\nu_1,\ldots,\nu_{\on{rk}\mathfrak{l}-\on{rk}\mathcal{D}(\mathfrak{l})}$ induce the isomorphism $L/\mathcal{D}(L) \iso (\BC^\times)^{\on{rk}\mathfrak{l}-\on{rk}\mathcal{D}(\mathfrak{l})}$.
\prop{}\label{gen_fund_gen}
$(a)$ Algebra $\CO(L)^L$ is generated (over $\BC$) by the following elements:
\begin{equation*}
\chi_{V_{L,i}},\, \chi_{\BC_{\nu_j}}^{\pm 1} \in \CO(L)^L,
\end{equation*}
where $i = 1,\ldots,\on{rk}\mathcal{D}(\mathfrak{l}),\, j=1,\ldots,\on{rk}\mathfrak{g}-\on{rk}\mathcal{D}(\mathfrak{l})$.

$(b)$ The differentials $\{d(\chi_{V_{L,i}}), d(\chi_{\BC_{\nu_j}})\,|\, i = 1,\ldots,\on{rk}\mathcal{D}(\mathfrak{l}),\, j=1,\ldots,\on{rk}\mathfrak{g}-\on{rk}\mathcal{D}(\mathfrak{l})\}$ are linearly independent at every regular element $C \in L^{\mathrm{reg}}$.
\eprop
\prf
Part $(a)$ is easy. 

Part $(b)$ follows from part $(a)$ together with Proposition \ref{gen_lie_centr} (use that $\on{dim}T_C Z_L(C)=\on{rk}\mathfrak{l}$).
\epr

\begin{remark}{}
{\em{Collection of representations $V_{L,1},\ldots, V_{L,\on{rk}\mathcal{D}(\mathfrak{l})}, \BC_{\nu_1},\ldots, \BC_{\nu_{\on{rk}\mathfrak{l}-\on{rk}\mathcal{D}(\mathfrak{l})}}$ as above can be though of as ``fundamental'' representations of $L$.}}
\end{remark}

Let now $G$ be a reductive group over complex numbers. Assume that $\mathcal{D}(G)=[G,G]$ is simply connected. Recall that $T \subset G$ is a maximal torus. Let $\al_{i}^{\vee},\, i \in I=\{1,\ldots,\on{rk}\mathcal{D}(\mathfrak{g})\}$ be  a set of simple coroots of $G$. 

Let $V_{1},\ldots, V_{\on{rk}\mathcal{D}(\mathfrak{g})}$ be irreducible (finite dimensional) modules over $G$. For $i=1,\ldots,\on{rk}\mathcal{D}(\mathfrak{g})$ let $\la_i^- \colon T \ra \BC^\times$ be the lowest weight of $V_{i}$. We assume that $\langle \la_i^-, \al_j^\vee \rangle = -\delta_{ij}$ i.e. that representations $V_i$  restrict to {\em{fundamental}} representations of $\mathcal{D}(G)$. Let $\BC_{\nu_1},\ldots,\BC_{\nu_{\on{rk}\mathfrak{g}-\on{rk}\mathcal{D}(\mathfrak{g})}}$ be one-dimensional representations of $G$ with weights $\nu_1,\ldots,\nu_{\on{rk}\mathfrak{g}-\on{rk}\mathcal{D}(\mathfrak{g})}$ such  that the characters $\nu_1,\ldots,\nu_{\on{rk}\mathfrak{g}-\on{rk}\mathcal{D}(\mathfrak{g})}$ induce the isomorphism $G/\mathcal{D}(G) \iso (\BC^\times)^{\on{rk}\mathfrak{l}-\on{rk}\mathcal{D}(\mathfrak{l})}$. Let $L \subset G$ be a (standard) Levi subgroup of $G$ that contains $T$. Let $J \subset I$ be the subset of simple roots corresponding to $L$. Let $V_{L,\la_i} \subset V_{\la_i}$ be the irreducible $L$-subrepresentation generated by the lowest vector of $V_{\la_i}$. 

\lem{}\label{restr_fund_is_fund}
Representations $V_{L,i}$ for $i \in J$ restrict to fundamental representations of $L$. Representations $\{V_{L,i},\,\BC_{\nu_{j}}\,|\, i \in I \setminus J,\, j=1,\ldots,\on{rk}\mathfrak{g}-\on{rk}\mathcal{D}(\mathfrak{g})\}$  are one-dimensional. The characters 
\begin{equation*}
\{\la^{-}_i,\, i \in I \setminus J\} \cup \{\nu_1,\ldots,\nu_{\on{rk}\mathfrak{g}-\on{rk}\mathcal{D}(\mathfrak{g})}\}
\end{equation*}
induce the isomorphism $L/\mathcal{D}(L) \iso (\BC^\times)^{\on{rk}\mathfrak{l}-\on{rk}\mathcal{D}(\mathfrak{l})}$. 
\elem
\prf
It is clear that  $V_{L,i}$ for $i \in J$ restrict to fundamental representations of $L$.  

Consider the exact sequence 
\begin{equation*}
1 \ra \mathcal{D}(L) \cap T \ra T \ra L/\mathcal{D}(L) \ra 1.
\end{equation*}
Restrictions of $\la^{-}_{i},\, i \in J$ to $\mathcal{D}(L) \cap T$ induce the isomorphism $\mathcal{D}(L) \cap T \iso (\BC^\times)^{\on{rk}\mathfrak{l}}$. It also follows from the definitions that the characters $\la^{-}_i$, $\nu_j$ induce the isomorphism $T \iso (\BC^\times)^r$. We conclude that $\la^{-}_i,\, i \in I \setminus J$, $\nu_1,\ldots,\nu_{\on{rk}\mathfrak{g}-\on{rk}\mathcal{D}(\mathfrak{g})}$ induce the isomorphism $L/\mathcal{D}(L) \iso (\BC^\times)^{\on{rk}\mathfrak{l}-\on{rk}\mathcal{D}(\mathfrak{l})}$.
\epr

\rem{}
{\em{Lemma \ref{restr_fund_is_fund} shows that ``fundamental'' weights of $G$ restrict to ``fundamental'' weights of $L$.}}
\erem


\prop{}\label{prop_lin_ind_inv}
The differentials  
\begin{equation*}
\{d(\chi_{V_{L,i}}),\, i=1,\ldots,\on{rk}\mathcal{D}(\mathfrak{g})\} \cup \{d(\chi_{\BC_{\nu_j}}),\, j=1,\ldots,\on{rk}\mathfrak{g}-\on{rk}\mathcal{D}(\mathfrak{g})\}
\end{equation*}
are linearly independent at every regular element $C \in L^{\mathrm{reg}}$. 
\eprop
\prf
Follows from Proposition \ref{gen_fund_gen} together with Lemma \ref{restr_fund_is_fund}.
\epr

\end{document}